\documentclass[11pt,a4paper]{article}

\usepackage{mathtools}
\usepackage{authblk} % author
\usepackage{algpseudocode,algorithmicx,algorithm}

\usepackage{mathrsfs}
\usepackage{latexsym,bm}
\usepackage{amsmath,amsfonts,amssymb,amsthm,rotating}
\usepackage{extarrows}

\usepackage{graphicx,subfigure,epstopdf,float}
\usepackage{enumerate,cases,multirow}
\usepackage{makecell}
\usepackage{caption}

\usepackage{longtable,colortbl,arydshln,threeparttable}
\definecolor{mygray}{gray}{.9}

\usepackage{indentfirst}
\usepackage[top=25mm,bottom=20mm,left=25mm,right=20mm]{geometry}
\usepackage{cite}

\usepackage{listings}

\usepackage{makeidx}        % generate an index, automatically
\usepackage{booktabs}
\usepackage{xcolor}
\usepackage[bookmarksnumbered,colorlinks=true,citecolor=red,linkcolor=red,hyperindex,linktocpage=true]{hyperref}

\newcommand{\ket}[1]{| #1 \rangle} % |u>
\newcommand{\bra}[1]{\langle #1 |} % <u|

\newcommand{\bb}{\boldsymbol}

\def \d {\mathrm{d}}
\def \e {\mathrm{e}}
\def \i {\mathrm{i}}

\newcounter{parentalgorithm}

\makeatother

\newtheorem{theorem}{Theorem}[section]
\newtheorem{lemma}{Lemma}[section]

\newtheorem{definition}{Definition}[section]

\theoremstyle{remark}
\newtheorem{remark}{\bf Remark}[section]

\numberwithin{equation}{section}

\begin{document}

\title{Schr\"odingerization based quantum algorithms for the fractional Poisson  equation}

\author[1]{Shi Jin\thanks{shijin-m@sjtu.edu.cn}}
\author[1, 2]{Nana Liu\thanks{nana.liu@quantumlah.org}}
\author[3]{Yue Yu\thanks{terenceyuyue@xtu.edu.cn}}
\affil[1]{School of Mathematical Sciences, Institute of Natural Sciences, MOE-LSC, Shanghai Jiao Tong University, Shanghai, 200240, China}
%\affil[2]{Shanghai Artificial Intelligence Laboratory, Shanghai, China}
\affil[2]{University of Michigan-Shanghai Jiao Tong University Joint Institute, Shanghai 200240, China}
\affil[3]{School of Mathematics and Computational Science, Hunan Key Laboratory for Computation and Simulation in Science and Engineering, Key Laboratory of Intelligent Computing and Information Processing of Ministry of Education, National Center for Applied Mathematics in Hunan, Xiangtan University, Xiangtan, Hunan 411105, China}

\maketitle

\begin{abstract}
  We develop a quantum algorithm for solving  high-dimensional fractional Poisson equations. By applying the Caffarelli-Silvestre extension, the \(d\)-dimensional fractional equation is reformulated as a local partial differential equation in \(d+1\) dimensions. We propose a quantum algorithm for the  finite element discretization of this local problem, by capturing the steady-state of the corresponding differential equations using the Schr\"odingerization approach from \cite{JLY22SchrShort, JLY22SchrLong, analogPDE}. The Schr\"odingerization technique transforms general linear partial and ordinary differential equations into Schr\"odinger-type systems, making them suitable for quantum simulation. This is achieved through the warped phase transformation, which maps the equation into a higher-dimensional space. We provide detailed implementations of the method and conduct a comprehensive complexity analysis, which can show up to exponential advantage -- with respect to the
inverse of the mesh size in high dimensions -- compared to its classical counterpart.
Specifically, while the classical method requires \(\widetilde{\mathcal{O}}(d^{1/2} 3^{3d/2} h^{-d-2})\) operations, the quantum counterpart requires \(\widetilde{\mathcal{O}}(d 3^{3d/2} h^{-2.5})\) queries to the block-encoding input models, with the quantum complexity being independent of the dimension \(d\) in terms of the inverse mesh size \(h^{-1}\).  Numerical experiments are conducted to verify the validity of our formulation.
\end{abstract}


%\tableofcontents

\section{Introduction}

The fractional Laplacian \( (-\Delta)^s \), with fractional order \( s \), is an important nonlocal generalization of the Laplacian operator \( -\Delta \). It has garnered increasing attention over the past few decades due to its ability to capture complex physical phenomena involving long-range interactions. This operator appears in various areas of application, such as anomalous diffusion, L\'evy processes, stochastic dynamics, image processing, groundwater solute transport, finance, and turbulent flows \cite{Klafter2005Anomalous,Carreras2001Anomalous,Gatto2014Fractional,Vazquez2014Fractional,Du2019Nonlocal,Caffarelli2007Fractional,Nochetto2015Fractional}. Since obtaining analytical solutions for differential equations involving the fractional Laplacian is challenging, the development of numerical methods has become essential. Various methods have been introduced to discretize problems with fractional Laplacians, including finite element methods (FEMs), spectral methods, Monte Carlo methods, meshless methods, finite difference methods, and deep learning approaches \cite{Acosta2017Fractional,Acosta2017codeFractional,Yang2023fastNonlocal,Hao2020Spectral,Sheng2023MC,Hao2025Meshless,
Hao2021fractionalFDM,YangChen2023fastNonlocalFDM,Gu2022DeepRitz}.

In this paper, we focus on quantum computation for this nonlocal problem, particularly in high-dimensional settings, an area that has not been explored to the best of our knowledge.
Let $\Omega$ be an open, bounded and connected subset of $\mathbb{R}^d$ ($d\ge 1$) with Lipschitz boundary $\partial \Omega$. Given $s\in (0,1)$ and a smooth enough function $f$, the fractional Poisson equation for $u$ is \cite{Caffarelli2007Fractional,Nochetto2015Fractional}
\begin{equation} \label{ProbFractional}
\begin{cases}
(-\Delta )^s u = f \qquad & \mbox{in}~~\Omega, \\
u = 0 \qquad & \mbox{on}~~\partial\Omega,
\end{cases}
\end{equation}
where the fractional Laplacian of a function $u: \mathbb{R}^d \to \mathbb{R}$ is defined by
\[(-\Delta )^s u = C_{d,s} \text{p.v.} \int_{\mathbb{R}^n} \frac{u(x) - u(\xi)}{|x-\xi|^{d + 2s}} \d \xi, \]
with $C_{d,s} = 4^s s \Gamma(s+d/2)/(\pi^{d/2}\Gamma(1-s))$ being the normalization constant and p.v. stands for the Cauchy principal value.

Due to the nonlocality and strong singularity, classical numerical simulations of fractional partial differential equations (PDEs) encounter more significant challenges from the curse of dimensionality compared to local PDEs in high dimensions.
Quantum computing offers a promising approach to overcome the computational challenges associated with high dimensionality, as it requires only a logarithmic number of qubits relative to the matrix dimension to store both the matrix and the solution vector. This advantage has led to increasing interest in developing quantum algorithms for solving PDEs \cite{Cao2013Poisson,Berry2014Highorder,qFEM-2016,Costa2019Wave,
Engel2019qVlasov,Childs2020qSpectral,Linden2020heat,Childs2021high,JinLiu2022nonlinear,GJL2022QuantumUQ,JLY2022multiscale}. Among the notable quantum strategies for solving PDEs is the {\it Schr\"odingerization} method introduced in \cite{JLY22SchrShort, JLY22SchrLong, analogPDE}. This simple and generic framework enables quantum simulation for {\it all} linear PDEs and ordinary differential equations (ODEs). The core idea is to apply a warped phase transformation that maps the equations to one higher dimension, which, in Fourier space, transforms them into a system of Schr\"odinger-type equations!
The approach has been expanded to address a wide array of problems, including open quantum systems within bounded domains with non-unitary artificial boundary conditions \cite{JLLY23ABC}, problems entailing physical boundary or interface conditions \cite{JLLY2024boundary,JLY24Circuits}, Maxwell's equations \cite{JLC23Maxwell,MJL2024time}, the Fokker-Planck equation \cite{JLY24FokkerPlanck}, multiscale PDEs \cite{HJZ2024multiscale}, stochastic differential equations \cite{JLW2024Stochastic}, ill-posed problems such as the backward heat equation \cite{JLM24SchrBackward}, linear dynamical systems with inhomogeneous terms \cite{JLM24SchrInhom}, non-autonomous ordinary and partial differential equations \cite{CJL23TimeSchr,CJL2024TimeSchr}, iterative linear algebra solvers \cite{JinLiu-LA}, etc. This approach is also natural for continuous variables quantum computing \cite{analogPDE} thus is also an idea~--~so far the only possible one~--~for the analog quantum simulation of PDEs and ODEs that completely retains the continuous nature of the ODEs or PDES without discretizations. This method can also be extended to parabolic partial differential equations using a Jaynes-Cummings-like model -- which is more widely experimentally accessible  -- as was proposed in \cite{JL24JaynesCummings}.

As demonstrated in the Schr\"odingerization approach, dimension lifting has proved to be a generic method that can be applied in the design of quantum algorithms for linear differential equations. This view point is also applied in various other contexts, such as transforming nonlinear PDEs into linear ones \cite{JinLiu2022nonlinear}, converting non-autonomous Hamiltonian systems with time-dependent Hamiltonians into autonomous PDEs with time-independent coefficients \cite{CJL23TimeSchr,CJL2024TimeSchr,JLY24Circuits,MJL2024time}, and developing new transformations that map PDEs with uncertainty to deterministic PDEs \cite{GJL2022QuantumUQ}.
While dimension lifting may present a challenge in classical computation, it does not pose the same issue in quantum computing thanks to the use of qubits.

In the present work, we employ the Caffarelli-Silvestre extension technique in \cite{Caffarelli2007Fractional} to the numerical solution of the fractional Poisson equation, which transforms it into a local elliptic PDE in one higher dimension, and subsequently develop a new quantum algorithm to solve it.
The algorithm involves three steps: 1) first, we apply the finite element method, along the lines of \cite{Nochetto2015Fractional}, to the corresponding local problem and obtain a linear
algebraic system $A\bb{x} = \bb{b}$; 2) second, we consider $\bb{x}$ as the steady state solution to a linear differential equation
$\frac{{\rm d} \bb{u}(t)}{{\rm d} t}  =	- A \bb{u}(t) + \bb{b}$; 3) finally, we calculate $\bb{u}(t)$ using the quantum simulation of ODEs via Schr\"odingerization.
A detailed complexity analysis reveals that the quantum algorithm provides an exponential advantage over its classical counterparts, especially in high-dimensional settings. Specifically, the quantum algorithm has a time complexity of \( \widetilde{\mathcal{O}}(d3^{\frac{3d}{2}}h^{-2.5}) \) for block-encoding input models, where $h$ is the spatial mesh size, while the classical Conjugate Gradient (CG) method scales as \( \widetilde{\mathcal{O}}(d^{\frac{1}{2}}3^{\frac{3d}{2}}h^{-d-2}) \), demonstrating  the exponential quantum advantage on $d$.

The paper is structured as follows. In Section \ref{sec:SchrAxb}, we introduce the Schr\"odingerization-based quantum subroutine to solve the linear systems problems and present the query complexity in terms of block-encoding input models. Section \ref{sec:CaffarelliSilvestre} provides an overview of the Caffarelli-Silvestre extension technique for the fractional Laplacian problem, as well as a review of standard finite element methods for the resulting (local) Poisson equation and an estimation of the matrix condition number. In Section \ref{sec:complexity}, we present the block encoding of the stiffness matrix on a hypercube domain, explain the process of extracting the solution vector, and analyze the query complexity for our Schr\"odingerization based approach. Numerical simulations are presented in Section \ref{sec:numerical}, where a nonlocal variational formulation of the original problem is also tested via the Schr\"odingerization approach. Conclusions are presented in the final section.

\section{Quantum algorithms} \label{sec:SchrAxb}

We will use linear finite element methods for numerically solving the fractional Poisson equation.  The specifics of the finite element formulations are provided in the next section. Once this is done, the remaining task is to solve the linear system $ A \bb{x} = \bb{b} $.

\subsection{The Schr\"odingerization based algorithm}

In the context of the problem we are considering, $ A $ is a symmetric positive definite matrix. In \cite{HJZ2024multiscale, JLMY2025SchOptimal}, the solution $ \bb{x} $ was interpreted as the steady state of the following linear ODEs:
\begin{equation}\label{reformulationODE}
	\frac{{\rm d} \bb{u}(t)}{{\rm d} t}  =
	- A \bb{u}(t) + \bb{b},\quad
	\bb{u} (0) = \bb{u}_0.
\end{equation}
This reformulation enables the application of the Schr\"odingerization to construct a Hamiltonian system for quantum computing. For simplicity, we choose $\bb{u}_0 = \bb{0}$ as the initial data, since the choice of the initial value does not affect the steady state.

%{\color{blue} The time complexity of this ODE-based method will have a quadratic dependence on the condition number $\kappa(A)$, regardless of whether the Schr\"odingerization approach or the LCHS method in \cite{ALL2023LCH,ACL2023LCH2} is used. } We can reduce the $\kappa$-dependence from quadratic to nearly linear by using the variable-time amplitude amplification (VTAA) in \cite{Childs2017QLSA}. However, VTAA necessitates substantial modifications to the LCU or QSP algorithms and introduces significant overhead, and the performance of VTAA in solving the quantum linear systems problems has not been quantitatively reported in the literature.

%We first demonstrate that directly applying Schr\"odingerization to equation \eqref{reformulationODE} results in a time complexity with quadratic dependence on the condition number $\kappa(A)$.

The solution to \eqref{reformulationODE} with $\bb{u}_0 = \bb{0}$ can be expressed as
\[\bb{u}(T) = \int_0^T \e^{-A (T-s)} \d s \bb{b} = \int_0^T \e^{-A s}\d s \bb{b} = (I - \e^{-A T}) A^{-1}\bb{b}=:S(A,T)\bb{b}, \]
where $T$ is the evolution time and $I$ is an identity matrix. One obtains from
\begin{equation}\label{approxOperator}
\|A^{-1} - S(A,T)\| = \| \e^{-A T}  A^{-1}\| = \frac{1}{|\lambda|_{\min}} \e^{-|\lambda|_{\min} T}
= \frac{\kappa}{\|A\|} \e^{-|\lambda|_{\min} T} \le \delta,
\end{equation}
which yields
\begin{equation}\label{Timesteady}
T \ge \frac{\kappa(A)}{\|A\|} \log \frac{\kappa(A)}{\|A\|\delta}.
\end{equation}
To avoid the complicated Linear Combination of Unitaries (LCU) procedure for the integration, we employ the augmentation technique in \cite{JLY22SchrLong,JLM24SchrInhom,JLMY2025SchOptimal} to derive a homogeneous ODE system by introducing a time-independent auxiliary vector:
\begin{equation}\label{reformulationODEaug}
	\frac{{\rm d}}{{\rm d} t} \bb{u}_f = A_f \bb{u}_f,
	\quad
	A_f = \begin{bmatrix}
		-A  & \frac{I}{T}\\
		O & O
	\end{bmatrix},\quad
	\bb{u}_f(0) =  \begin{bmatrix}
		\bb{0} \\
		T \bb{b}
	\end{bmatrix}.
\end{equation}
It is straightforward to verify that \eqref{reformulationODEaug}, with solution given by
$\bb{u}_f = [\bb{u}^\top, T \bb{b}^\top]^\top$, is equivalent to \eqref{reformulationODE}, where ``$\top$'' stands for the matrix transpose.

For the Schr\"odingerization of the ODE system, we first decompose $A_f$ into a Hermitian term and an anti-Hermitian term $A_f = H_1 + \i H_2$, with
\begin{equation*}
	H_1 = \frac{1}{2}(A_f + A_f^{\dagger}) = \begin{bmatrix}
		-A& \frac{I}{2T} \\
		\frac{I}{2T}  & O
	\end{bmatrix},\quad
	H_2 = \frac{1}{2\i}(A_f - A_f^{\dagger}) =  \frac{1}{2\i}
	\begin{bmatrix}
		O & \frac{I}{T} \\
		-\frac{I}{T} &O
	\end{bmatrix},
\end{equation*}	
where ``$\dagger$'' is the conjugate transpose, and both $H_1$ and $H_2$ are Hemitian.
Using the warped phase transformation $\bb{v}(t,p) = \e^{-p}\bb{u}_f(t)$ for $p\ge 0$, with the equation extended naturally to $p<0$, we get
\begin{equation}\label{uf2v}
	\frac{{\rm d}}{{\rm d} t} \bb{v}(t) =  -H_1 \partial_p \bb{v} + \i H_2 \bb{v}.
\end{equation}
The initial data for \eqref{uf2v} is specified as $\bb{v}(0,p) = \zeta(p) \bb{u}_f(0)$, where $\zeta(p) \in H^r(\mathbb{R})$, with $r \geq 1$, is chosen under the condition the condition that $\zeta(p) = \e^{-p}$ for $p \geq 0$ and decays sufficiently fast as $p\to -\infty$.  A particular case is $\zeta(p) = \e^{-|p|} \in H^1(\mathbb{R})$. A construction of $\zeta$ with higher regularity can be found in \cite{JLM24SchrBackward}. When $A$ is positive semi-definite, $\bb{u}_f(t)=\e^p \bb{v}(t,p)$ for all $p>p_* = \frac12$ (see Theorem 3.5 in \cite{JLMY2025SchOptimal}),  so we can recover \(\bb{u}_f(t)\) by projecting $\bb{v}$ onto \(p > p_*\), and then obtain the solution \(\bb{u}(t)\) by extracting the first block of \(\bb{u}_f(t)\).

Let $p\in [-L,R]$ with $L, R>0$ and $R>0$ satisfying $\e^{-L} \approx 0$ and $\e^{-R} \approx 0$. Then one can apply the periodic boundary condition in the $p$ direction and use the Fourier spectral method by discretizing the $p$ domain. The choice of $L$ and $R$ was provided in Theorem 3.5 of \cite{JLMY2025SchOptimal}, given by
\[L =\mathcal{O}(\kappa(A)\log \frac{1}{\varepsilon}+\frac 1 2),\quad  R = \mathcal{O} \Big(\log \frac{1}{\varepsilon} +\frac{1}{2}\Big),\]
where $\varepsilon$ is the desired precision. We choose a uniform mesh size $\Delta p = (L+R)/N_p$ for the auxiliary variable with $N_p=2^{n_p}$ being an even number. The grid points are denoted by $-L = p_0<p_1<\cdots<p_{N_p} = R$. For $p\in [-L, R]$, the one-dimensional basis functions for the Fourier spectral method are usually chosen as
\[\phi_l(p) = \e ^{\i \mu_l (p+L)} , \quad \mu_l = \frac{2\pi (l-N_p/2) }{R+L}, \quad l=0,1,\cdots,N_p-1.\]
We also define
\[\Phi = (\phi_{jl})_{N_p\times N_p} = (\phi_l(p_j))_{N_p\times N_p}, \qquad D_\mu = \text{diag} ( \mu_0, \cdots, \mu_{N_p-1} ).\]

Let $\bb v(t,p) = [v_1(t,p), \cdots, v_N(t,p)]^\top$ be the solution to \eqref{uf2v}. The spectral discretization of $v_i(t,p)$ is
\[v_{i,h}(t,p) = \sum_{l=0}^{N_p-1} \tilde{v}_{i,l,h}(t) \phi_l(p), \qquad i = 1,\cdots,N,\]
where we use the subscript $h$ to denote the numerical solution for \eqref{uf2v}. The approximate solution is then given by $\bb{v}_h(t,p) = [v_{1,h}(t,p),\cdots, v_{N,h}(t,p)]^\top$, which can be written as
\begin{equation}\label{interpcoeff}
\bb{v}_h(t,p) = \sum_{l=0}^{N_p-1} \tilde{\bb{v}}_{l,h}(t) \phi_l(p), \qquad
	\tilde{\bb{v}}_{l,h}(t) = \frac{1}{N_p} \sum\limits_{k=0}^{N_p-1} \bb{v}_h(t,p_k) \e ^{ - \i \mu_l (p_k+L)} .
\end{equation}
Let the vector $\bb{w}_h$ be the collection of the function $\bb v_h$ at the grid points, defined more precisely as
\[\bb{w}_h(t) = \sum_{k,i} v_{i,h}(t,p_k) \ket{k,i}
= [\bb{v}_h(t,p_0) ; \cdots; \bb{v}_h(t,p_{N_p-1})].\]
Accordingly, we define
\[\tilde{\bb{w}}_h(t) = \sum_{l,i} \tilde{v}_{i,l,h}(t) \ket{l,i}
= [\tilde{\bb{v}}_{0,h}(t) ; \cdots; \tilde{\bb{v}}_{N_p-1,h}(t)].\]
One can find that $\bb{w}_h(t) = (\Phi \otimes I_{N\times N}) \tilde{\bb{w}}_h(t)$. The ODEs are then transferred to a Hamiltonian system
\begin{equation} \label{eq:hamiltonian}
\frac{\rm d}{{\rm d} t} \bb{w}_h = -\i(P_\mu \otimes H_1) \bb{w}_h + \i (I_{n_p} \otimes H_2)\bb{w}_h, \qquad
\bb{w}_h(0) = \bb{\zeta} \otimes \bb{u}_f(0),
\end{equation}
where $I_{n_p}$ is an $n_p$-qubit identity matrix and $\bb{\zeta} = [\zeta(p_0), \cdots, \zeta(p_{N_p-1})]^\top$.
In terms of $\tilde{\bb{w}}_h = (\Phi^{-1} \otimes I)\bb{w}_h$, one gets
\begin{equation}\label{discreteLCHSQLSP}
\begin{cases}
\dfrac{\d}{\d t} \tilde{\bb{w}}_h(t) = -\i (D_\mu \otimes H_1 - I_{n_p} \otimes H_2) \tilde{\bb{w}}_h(t), \\
\tilde{\bb{w}}_h(0) = \tilde{\bb{\zeta}} \otimes \bb{u}_f(0), \quad \tilde{\bb{\zeta}} = \Phi^{-1}\bb{\zeta}.
\end{cases}
\end{equation}

\subsection{The complexity of Schr\"odingerization based algorithm}

Block encoding is a general input model for matrix operations on a quantum computer \cite{Gilyen2019QSVD,Chakraborty2019blockEncode,Lin2022Notes,ACL2023LCH2,JLMY2025SchOptimal}. Let $ A $ be an $ n $-qubit matrix. While $ A $ is not necessarily a unitary operator, there exists a unitary matrix $ U_A $ on $ (n+m) $-qubits such that the matrix $ \bar{A} = A / \alpha $ is encoded in the upper-left block of $ U_A $, where $ \|A\| \le \alpha $.
The effect of $ U_A $ acting on a quantum state can be described as
\[
U_A \ket{0^m, b} = \ket{0^m} \bar{A} \ket{b} + \ket{\bot},
\]
where $ \ket{0^m} = \ket{0}_m = \ket{0} \otimes \cdots \otimes \ket{0} $ is the $ m $-qubit state and $ \ket{\bot} $ is a component orthogonal to $ \ket{0^m} \bar{A} \ket{b} $. In practice, it is sufficient to find $U_A$ to block encode $A$ up to some error $\varepsilon$.

\begin{definition}\label{def:blockencoding}
	Suppose that $ A $ is an $ n $-qubit matrix and let $ \Pi = \bra{0^m} \otimes I $ with $ I $ being an $ n $-qubit identity matrix. If there exist positive numbers $ \alpha $ and $ \varepsilon $, as well as a unitary matrix $ U_A $ of $ (m+n) $-qubits, such that
	\[
	\| A - \alpha \Pi U_A \Pi^\dag \| = \| A - \alpha (\bra{0^m} \otimes I) U_A (\ket{0^m} \otimes I) \| \le \varepsilon,
	\]
	then $ U_A $ is called an $ (\alpha, m, \varepsilon) $-block-encoding of $ A $.
\end{definition}

For a sparse matrix, we can construct the block encoding using its sparse query model \cite{Gilyen2019QSVD, Chakraborty2019blockEncode, Lin2022Notes}, with the details omitted.
Compared to sparse encoding, block-encoding serves not only as an input model for quantum algorithms but also enables various matrix operations to implement block-encoding of more complex matrices. For instance, with the coefficient matrix derived from tensor-product finite elements in Section \ref{subsec:BE}, we can construct the block-encoding for high-dimensional problems by combining block-encodings of one-dimensional coefficient matrices, thereby avoiding explicit computation of non-zero elements in high-dimensional coefficient matrices. This approach is particularly advantageous since the sparsity of high-dimensional coefficient matrices typically grows exponentially with dimension. Furthermore, block-encoding remains applicable to dense matrices where sparse encoding fails~-~as seen in the nonlocal variational formulation discussed later (see Section \ref{subsec:nonlocal}). However, as noted in Remark \ref {rem:fast}, for regular domains the nonlocal form may still admit efficient block-encoding through specialized techniques.

Following the similar implementation in \cite{JLMY2025SchOptimal} with some modifications, we may conclude the following lemma.

\begin{lemma}\label{lem:Schrcomplexity}
Suppose that the coefficient matrix of  \eqref{reformulationODE} is negative semi-definite over the interval $[0,T]$. Let $\|A\| \le \alpha_A$. Then there exists a quantum algorithm that prepares an $\epsilon$-approximation of the state $\ket{\bb{u}(T)}$ with $\Omega(1)$ success probability and a flag indicating success, using
\[\widetilde{\mathcal{O}}\Big(\frac{\|\bb{u}(0)\| + T \|\bb{b}\|}{\|\bb{u}(T)\|} \alpha_A \eta_{\max} T \log(1/\varepsilon)\Big)\]
queries to the block-encoding oracle of $A$ and
\[\mathcal{O}\Big(\frac{\|\bb{u}(0)\| + T \|\bb{b}\|}{\|\bb{u}(T)\|} \Big)\]
queries to the state preparation oracles for $\tilde{\bb{w}}_h(0)$ and $\bb{b}$. Here, $\eta_{\max} = \max_k |\eta_k|$ represents the maximum absolute value among the discrete Fourier modes, which is $\mathcal{O}(\log (1/\varepsilon))$ if we choose a sufficiently smooth initialization of $\zeta(p)$.
\end{lemma}

\begin{theorem}\label{thm:SchrCost}
Let $A$ be a positive definite matrix and $\alpha_A = \Omega(\|A\|)$. If we choose
\[T = \Theta\Big(\frac{\kappa}{\|A\|}  \sqrt{2 \log\frac{4\kappa}{\xi \|A\|\varepsilon}}\Big), \qquad \xi = \|A^{-1} \ket{\bb{b}}\|,\]
then there exists a quantum algorithm that prepares an $\epsilon$-approximation of the state $\ket{\bb{x}}$ with $\Omega(1)$ success probability and a flag indicating success, using
\[\widetilde{\mathcal{O}}\Big( \frac{\kappa^2}{\xi \|A\|}  \log^2 \frac{\kappa}{\xi \|A\| \varepsilon}\log\frac{1}{\varepsilon}\Big)\]
queries to the block-encoding oracle of $A$ and
\[\mathcal{O}\Big( \frac{\kappa}{\xi \|A\|}  \log^{0.5}\frac{\kappa}{\xi \|A\| \varepsilon} \Big)\]
queries to the state preparation oracle for $\bb{b}$.
\end{theorem}
\begin{proof}
Let $\bb{x}_T = \bb{u}(T)$. By Lemma \ref{lem:Schrcomplexity}, $\|\ket{\bb{x}_T} - \ket{\bb{x}_T^{\text{d}}}\| \le \varepsilon$, so we can require $\|\ket{\bb{x}} - \ket{\bb{x}_T}\|\le \varepsilon$.
According to Eq.~\eqref{approxOperator},
\[\|\bb{x} - \bb{x}_T\|
= \|A^{-1}\bb{b} - S(A,T)\bb{b}\| \le \|A^{-1} - S(A,T)\|\|\bb{b}\| \le \delta \|\bb{b}\|.\]
Using the inequality $\| \frac{\bb{x}}{\|\bb{x}\|} - \frac{\bb{y}}{\|\bb{y}\|} \| \le 2 \frac{\|\bb{x} - \bb{y}\|}{\|\bb{x}\|}$ for two vectors $ \bb{x}, \bb{y} $, we can bound the error as
\begin{align*}
\|\ket{\bb{x}} - \ket{\bb{x}_T}\|
 \le 2 \frac{\|\bb{x} - \bb{x}_T\|}{\|\bb{x}\|} \le 2 \frac{\delta \|\bb{b}\|}{\|\bb{x}\|} \le \varepsilon,
\end{align*}
which gives
\[\delta = \frac{\|\bb{x}\|}{2\|\bb{b}\|} \varepsilon =  \frac{\|A^{-1} \ket{\bb{b}}\|}{2} \varepsilon = : \frac{\xi\varepsilon}{2}, \qquad \xi = \|A^{-1} \ket{\bb{b}}\|. \]

When $\bb{u}(0) = \bb{0}$, one has
\[\frac{\|\bb{u}(0)\| + T \|\bb{b}\|}{\|\bb{u}(T)\|}
\approx \frac{T \|\bb{b}\|}{\|\bb{x}\|} = \frac{T }{\|A^{-1}\ket{\bb{b}}\|} =\frac{T}{\xi}
= \mathcal{O}\Big( \frac{\kappa}{\xi \|A\|}  \log^{0.5}\frac{\kappa}{\xi \|A\| \varepsilon} \Big)\]
and
\[\frac{\|\bb{u}(0)\| + T \|\bb{b}\|}{\|\bb{u}(T)\|} \alpha_A \eta_{\max} T \log(1/\varepsilon)
= \mathcal{O}\Big( \frac{\kappa^2}{\xi \|A\|}  \log^2 \frac{\kappa}{\xi \|A\| \varepsilon}\eta_{\max} \Big).\]
This completes the proof.
\end{proof}

\begin{remark}
The factor $\xi \|A\| = \|A\| \|A^{-1} \ket{\bb{b}} \| \in [1, \kappa]$ is invariant under the scaling of $A$ (i.e, $A \to c A$).
\end{remark}

\begin{remark}\label{rem:VTAA}
We can reduce the dependence on $ \kappa $ from quadratic to nearly linear by using variable-time amplitude amplification (VTAA), as demonstrated in \cite{Childs2017QLSA}. However, VTAA requires substantial modifications to the LCU or Quantum Signal  Processing (QSP) algorithms and introduces significant overhead itself.
A simpler approach to addressing this issue within ODE-based framework for solving Quantum Linear System Problem (QLSP) remains an open question.
\end{remark}

\section{The Caffarelli-Silvestre extension for fractional Laplacian problem} \label{sec:CaffarelliSilvestre}

This section reviews a dimension-lifting technique for the numerical solution of the fractional Laplacian problem.

\subsection{The Caffarelli-Silvestre extension}

The fractional Laplacian is a nonlocal operator, which is one of the main difficulties in the study of problem \eqref{ProbFractional} in  both theoretical and numerical aspects. Caffarelli and Silvestre found in \cite{Caffarelli2007Fractional} that the fractional Laplacian in $\mathbb{R}^d$ can be realized as a Dirichlet-to-Neumann (DtN) map via an extension problem on the upper half-space $\mathbb{R}_+^{d+1}$. This leads to a ({\it local}) elliptic PDE in one higher dimension. Building on this, Nochetto et al. systematically examined the solution technique within a bounded domain and numerically solved the resulting problem using a first-degree tensor product finite element discretization in \cite{Nochetto2015Fractional}.

Before presenting the Caffarelli-Silvestre extension, we first introduce some notations used in \cite{Nochetto2015Fractional}.
It is well-known that the operator $-\Delta: L^2(\Omega) \to L^2(\Omega)$, with the domain $\text{Dom}(-\Delta) = \{v \in H_0^1(\Omega): \Delta v \in L^2(\Omega) \}$, is positive, unbounded, and closed, and its inverse is compact. This implies there exist $\{\lambda_k, \varphi_k\}_{k \in \mathbb{N}} \subset \mathbb{R}_+ \times H_0^1(\Omega)$ such that $\{\varphi_k\}_{k \in \mathbb{N}}$ is an orthogonal basis of $H_0^1(\Omega)$ and
\begin{equation}\label{eigenLap}
-\Delta \varphi_k = \lambda_k \varphi_k \quad \text{in} ~~ \Omega,  \qquad  k \in \mathbb{N},
\end{equation}
where $\lambda_1 <\lambda_2 < \cdots$.
The fractional powers of the Dirichlet Laplacian $(-\Delta)^s$ can be defined for $u \in C_0^{\infty}(\Omega)$ by
\[(-\Delta)^s u = \sum_{k = 1}^{\infty} u_k \lambda_k^s \varphi_k, \qquad u_k = \int_{\Omega} u \varphi_k \d x. \]
This operator can be extended to the Hilbert space
\[\mathbb{H}^s(\Omega) = \Big\{w = \sum_{k=1}^\infty w_k \varphi_k \in L^2(\Omega):  \|w\|_{\mathbb{H}^s(\Omega)}^2 = \sum_{k=1}^\infty \lambda_k^s|w_k|^2 < \infty \Big\},\]
which can be characterized by fractional Sobolev space as
\[\mathbb{H}^s(\Omega) = \begin{cases}
H^s(\Omega), \qquad & s \in (0,\frac12), \\
H_{00}^{1/2}(\Omega), \qquad &s = \frac12, \\
H_0^s(\Omega), \qquad &s \in (\frac12, 1).
\end{cases}\]
The dual space of $\mathbb{H}^s(\Omega)$ is denoted by $\mathbb{H}^{-s}(\Omega)$.
Here, for $0<s<1$, the Sobolev space $H^s(\Omega)$ is defined by
\[H^s(\Omega) = \{ w\in L^2(\Omega): |w|_{H^s(\Omega)} < \infty\},\]
equipped with the norm
\[\|w\|_{H^s(\Omega)} = (\|w\|_{L^2(\Omega)}^2 +   |w|_{H^s(\Omega)}^2 )^{1/2},\]
where
\[|w|_{H^s(\Omega)} := \int_{\Omega}\int_{\Omega} \frac{w(x) - w(x')}{|x-x'|^{d + 2s}} \d x  \d x'.\]
The space $H_0^s(\Omega)$ is defined as the closure of $C_0^{\infty}(\Omega)$ with respect to the norm $\|\cdot\|_{s,\Omega}$. Note that when $\partial \Omega$ is Lipschitz, $H_0^s(\Omega) = H^s(\Omega)$ if $s \le \frac12$ since $C_0^{\infty}(\Omega)$ is dense in $H^s(\Omega)$ if and only if $s \le \frac12$. The space $H_{00}^{1/2}(\Omega)$ is the so-called Lions–Magenes space and can be defined by an interpolation space, that is, $H_{00}^{1/2}(\Omega) = [H_0^1(\Omega), L^2(\Omega)]_{\frac12}$.

For the given domain $\Omega$, we define the semi-infinite cylinder
\[\mathcal{C} = \Omega \times (0,\infty)\]
and its lateral boundary
\[\partial_L\mathcal{C} = \partial\Omega \times [0,\infty).\]
Consider a function $u$ defined on $\Omega$. The $\alpha$-harmonic extension of $u(x)$ to the cylinder $\mathcal{C}$ is defined by the function $\mathfrak{u}(x,y)$ that solves the boundary value problem
\begin{equation}\label{alphaHarmonic}
\begin{cases}
\text{div}( y^\alpha \nabla \mathfrak{u}) = 0 \quad & \text{in}~~\mathcal{C}, \\
\mathfrak{u} = 0 \quad & \text{on}~~\partial_L \mathcal{C}, \\
\mathfrak{u} = u \quad & \text{on}~~\Omega \times \{0\}.
\end{cases}
\end{equation}

\begin{lemma}[The Caffarelli-Silvestre extension, \cite{Caffarelli2007Fractional,Nochetto2015Fractional}]
If $s\in (0,1)$ and $u \in \mathbb{H}^s(\Omega)$, then the $\alpha$-harmonic extension $\mathfrak{u}$ satisfies
\[d_s (-\Delta)^s u = \frac{\partial \mathfrak{u}}{\partial \nu^\alpha} = - \lim_{y \to 0^+} y^\alpha \mathfrak{u}_y\in \mathbb{H}^{-s}(\Omega)\]
in the sense of distributions, where $\alpha = 1-2s$ and
\[d_s = 2^{1-2s} \frac{\Gamma(1-s)}{\Gamma(s)}.\]
\end{lemma}

To exploit the Caffarelli-Silvestre extension, we need to address the degeneracy or singularity arising from $ y^\alpha$. To this end, we consider weighted Sobolev spaces with the weight function $y^\alpha$, where $\alpha \in (-1,1)$. Let $\mathcal{D} \subset \mathbb{R}^{d+1}$ be an open set. We define $w(x,y)\in L^2(\mathcal{D}, y^{\alpha})$  as a measurable function defined on $\mathcal{D}$ such that
\[\|w\|_{L^2(\mathcal{D}, y^{\alpha})}^2 = \int_{\mathcal{D}} |y|^{\alpha} w^2 \d x \d y < \infty.\]
Then we define the weighted Sobolev space
\[H^1(\mathcal{D}, y^{\alpha}) = \{w\in L^2(\mathcal{D}, y^{\alpha}): |\nabla w| \in L^2(\mathcal{D}, y^{\alpha}) \},\]
which, equipped with the norm
\[\|w\|_{H^1(\mathcal{D}, y^{\alpha})} = ( \|w\|_{L^2(\mathcal{D}, y^{\alpha})}^2 + \|\nabla w\|_{L^2(\mathcal{D}, y^{\alpha})}^2)^{1/2},\]
is a Hilbert space. Define
\[\mathring{H}_L^1(\mathcal{C}, y^{\alpha}) = \{w\in  H^1(\mathcal{C}, y^{\alpha}): w = 0 ~~ \text{on}~~\partial_L \mathcal{C}\}.\]

\begin{lemma}\cite{Nochetto2015Fractional}
Given $f\in \mathbb{H}^{-s}(\Omega)$, a function $u \in \mathbb{H}^s(\Omega)$ solves \eqref{ProbFractional} if and only if its $\alpha$-harmonic extension $\mathfrak{u}\in \mathring{H}_L^1(\mathcal{C}, y^{\alpha})$ solves
\begin{equation}\label{Probextension}
\begin{cases}
\text{div}( y^\alpha \nabla \mathfrak{u}) = 0 \quad & \text{in}~~\mathcal{C}, \\
\mathfrak{u} = 0 \quad & \text{on}~~\partial_L \mathcal{C}, \\
\frac{\partial \mathfrak{u}}{\partial \nu^\alpha} = d_s f \quad & \text{on}~~\Omega \times \{0\}.
\end{cases}
\end{equation}
This implies $u(x) = \mathfrak{u}(x,0)$ when $\mathfrak{u}$ solves \eqref{Probextension}.
\end{lemma}

The variational formulation for the extension problem \eqref{Probextension} is to find $\mathfrak{u}\in \mathring{H}_L^1(\mathcal{C}, y^{\alpha})$ such that
\begin{equation}\label{varExtension}
 \int_{\mathcal{C}} y^{\alpha} \nabla \mathfrak{u} \cdot \nabla \phi\, \d x \d y = d_s \langle f, \text{tr}\phi \rangle_{\mathbb{H}^{-s}(\Omega)\times \mathbb{H}^s(\Omega)}, \qquad \phi \in \mathring{H}_L^1(\mathcal{C}, y^{\alpha}),
\end{equation}
where $\langle \cdot, \cdot \rangle_{\mathbb{H}^{-s}(\Omega)\times \mathbb{H}^s(\Omega)}$ denotes the duality pairing between $\mathbb{H}^{-s}(\Omega)$ and $\mathbb{H}^s(\Omega)$. Here, for a function $w\in H^1(\mathcal{C}, y^{\alpha})$,  we denote by $\text{tr}w$ the trace of $w$ onto $\Omega \times \{0\}$. According to Proposition 2.5 of \cite{Nochetto2015Fractional},  $\text{tr}\mathring{H}_L^1(\mathcal{C}, y^{\alpha}) = \mathbb{H}^s(\Omega)$.

\subsection{Truncation of the extension problem}

For numerical implementation, we should truncate the semi-infinite cylinder $\mathcal{C}$ to a finite cylinder $\mathcal{C}_{\mathcal{Y}}$ for a suitably defined $\mathcal{Y}>0$, where $\mathcal{C}_{\mathcal{Y}} = \Omega \times (0, \mathcal{Y})$
with the lateral boundary $\partial _L \mathcal{C}_{\mathcal{Y}} = \partial \Omega \times [0, \mathcal{Y}]$.

As shown in \cite[Proposition 3.1]{Nochetto2015Fractional}, the solution $\mathfrak{u}$ decays exponentially in the extended direction.
%\begin{lemma}
%For every $\mathcal{Y}\ge 1$, the solution to \eqref{Probextension} satisfies
%\[\|\nabla \mathfrak{u}\|_{L^2(\Omega \times (\mathcal{Y}, \infty), y^\alpha)}
%\lesssim \e^{-\sqrt{\lambda_1} \mathcal{Y}/2} \|f\|_{\mathbb{H}^{-s}},\]
%where $\lambda_1$ is defined in \eqref{eigenLap}.
%\end{lemma}
This motivates the approximation of $\mathfrak{u}$ by a function ${\mathfrak{u}_{\mathcal{Y}}}$ that solves
\begin{equation}\label{Probextensiontruncated}
\begin{cases}
\text{div}( y^\alpha \nabla {\mathfrak{u}_{\mathcal{Y}}}) = 0 \quad & \text{in}~~\mathcal{C}_{\mathcal{Y}}, \\
{\mathfrak{u}_{\mathcal{Y}}} = 0 \quad & \text{on}~~\Gamma_D, \\
\frac{\partial {\mathfrak{u}_{\mathcal{Y}}}}{\partial \nu^\alpha} = d_s f \quad & \text{on}~~\Omega \times \{0\}
\end{cases}
\end{equation}
for a sufficiently large $\mathcal{Y}$, where
\begin{equation}\label{GammaD}
\Gamma_D = \partial_L \mathcal{C}_{\mathcal{Y}} \cup (\Omega \times \{ \mathcal{Y} \})
\end{equation}
 is called the the Dirichlet boundary.
For this truncated problem, we define
\[\mathring{H}_L^1(\mathcal{C}_{\mathcal{Y}}, y^{\alpha}) = \{w\in  H^1(\mathcal{C}, y^{\alpha}): w = 0 ~~ \text{on}~~\Gamma_D\}.\]
Note that the functions in $\mathring{H}_L^1(\mathcal{C}_{\mathcal{Y}}, y^{\alpha})$ can be viewed as elements of $\mathring{H}_L^1(\mathcal{C}, y^{\alpha})$ through extension by zero for $y >\mathcal{Y}$.
The variational problem is to find $\mathfrak{u}_{\mathcal{Y}}\in \mathring{H}_L^1(\mathcal{C}_{\mathcal{Y}}, y^{\alpha})$ such that
\begin{equation}\label{varExtensiontruncated}
 \int_{\mathcal{C}_{\mathcal{Y}}} y^{\alpha} \nabla {\mathfrak{u}_{\mathcal{Y}}} \cdot \nabla \phi \d x \d y = d_s \langle f, \text{tr} \phi \rangle_{\mathbb{H}^{-s}(\Omega)\times \mathbb{H}^s(\Omega)}, \qquad \phi \in \mathring{H}_L^1(\mathcal{C}_{\mathcal{Y}}, y^{\alpha}).
\end{equation}
The approximation property of ${\mathfrak{u}_{\mathcal{Y}}}$ is described below.

\begin{lemma}\cite[Theorem 3.5]{Nochetto2015Fractional}
For every $\mathcal{Y}\ge 1$, there holds
\[\|\nabla (\mathfrak{u}-{\mathfrak{u}_{\mathcal{Y}}})\|_{L^2(\mathcal{C}, y^\alpha)}
\le  C \e^{-\sqrt{\lambda_1} \mathcal{Y}/4} \|f\|_{\mathbb{H}^{-s}},\]
where $\lambda_1$ is defined in \eqref{eigenLap} and $C$ depends only on $s$ and $\Omega$. In particular, for every $\epsilon>0$, let
\[\mathcal{Y}_0 = \frac{2}{\sqrt{\lambda}_1} \Big(\frac12 \log C + 2 \log \frac{1}{\epsilon}\Big).\]
Then, for $\mathcal{Y} \ge \max\{ \mathcal{Y}_0, 1\}$, we have
\[\|\nabla (\mathfrak{u}-{\mathfrak{u}_{\mathcal{Y}}})\|_{L^2(\mathcal{C}, y^\alpha)}
\le  \epsilon \|f\|_{\mathbb{H}^{-s}}.\]
\end{lemma}

\subsection{Finite element discretization}

 Let $\mathcal{T}_{\mathcal{Y}}$ be a mesh of $\overline{\Omega} \times [0, \mathcal{Y}]$.
The generic element of $\mathcal{T}_{\mathcal{Y}}$ is denoted by $E = K \times I$, where $K \subset \mathbb{R}^d$ is an element isoparametrically equivalent to the unit cube $[0,1]^d$ and $I$ is an interval in the extended dimension. The finite element is selected as the tensor product element, with the space defined by
\[V_h = \{ w \in C^0(\overline{\mathcal{C}_{\mathcal{Y}}}): w|_E \in \mathbb{Q}_1(E),~~E\in \mathcal{T}_{\mathcal{Y}}, ~~w|_{\Gamma_D} = 0\},\]
where $\Gamma_D$ is defined in \eqref{GammaD}.
The discrete variational problem is: Find ${(\mathfrak{u}_{\mathcal{Y}})}_h \in V_h$ such that
\begin{equation}\label{FEMExtension}
 \int_{\mathcal{C}_{\mathcal{Y}}} y^{\alpha} \nabla {(\mathfrak{u}_{\mathcal{Y}})}_h \cdot \nabla w_h\, \d x \d y = d_s \langle f, \text{tr} w_h \rangle , \qquad w_h \in V_h.
\end{equation}
The finite element approximation of $u\in \mathbb{H}^s(\Omega)$ is then given by
\[u_h = \text{tr}{(\mathfrak{u}_{\mathcal{Y}})}_h.\]

For the extension problem \eqref{Probextensiontruncated}, the coefficient $ y^{\alpha} $, where $ \alpha = 1 - 2s $, either degenerates ($ s < 1/2 $) or blows up ($ s > 1/2 $). This necessitates the use of anisotropic elements to accurately capture the degenerate or singular behavior of the solution $ \mathfrak{u} $ at $ y \approx 0^+ $. In what follows, we consider the graded meshes used in \cite{Nochetto2015Fractional}. The graded partition for the interval $[0, \mathcal{Y}]$ is given by
\begin{equation}\label{graded}
y_k = \Big(\frac{k}{N} \Big)^\gamma \mathcal{Y}, \quad k = 0,\cdots,N,
\end{equation}
where $\gamma = \gamma(\alpha) > \frac{3}{1-\alpha} > 1$. In addition, we assume $\mathcal{T}_\Omega$ to be a shape-regular and quasi-uniform partition of $\Omega$, with $\sharp \mathcal{T}_\Omega \approx N^d$. The mesh $\mathcal{T}_{\mathcal{Y}}$ is constructed by the tensor product of $\mathcal{T}_\Omega$ and the graded partition. This implies $\sharp \mathcal{T}_{\mathcal{Y}} \approx N^{d+1}$ and $h_{\Omega} \approx N^{-1}$, where $h_\Omega$ is the mesh size of $\mathcal{T}_\Omega$.

For the graded meshes, we have the following error estimates in terms of $u$.
\begin{lemma}\label{lem:err}
Let $u$ and $u_h$ be the exact and numerical solutions to the fractional Poisson  problem \eqref{ProbFractional}. If $f\in \mathbb{H}^{1-s}$, which implies $u \in \mathbb{H}^{1+s}$, then
\[\|u - u_h\|_{\mathbb{H}^s(\Omega)} \lesssim | \log (\sharp \mathcal{T}_{\Omega}) |^s (\sharp \mathcal{T}_{\Omega})^{-1/d} \|u\|_{\mathbb{H}^{1+s}(\Omega)}
\approx | d \log (h_{\Omega}^{-1} )| ^s h_{\Omega} \|u\|_{\mathbb{H}^{1+s}(\Omega)}. \]
\end{lemma}

We now present the resulting linear system. For simplicity, we assume that $\Omega = (0,1)^d$ in what follows. The basis of $V_h$ is given by
\[\phi_{\bm{i}}(x,y) = \varphi_{i_1}(x_1) \cdots \varphi_{i_d}(x_d) \psi_{i_{d+1} }(y) ,\]
where $\bm{i} = (i_1,\cdots,i_d, i_{d+1})$, $\varphi_{i_k}$ is the piecewise linear function along $x_k$-direction for $k=1,\cdots,d,d+1$ and $\psi_{i_{d+1} }(y)$ is the piecewise linear function in the extended dimension. Then we can derive a linear system
\begin{equation}\label{AxbFEM}
A \bm{u}_{\mathcal{Y}} = \bm{b},
\end{equation}
where the stiffness matrix $A$ and force vector $\bb{b}$ are defined by
\[A = (A_{\bm{i}, \bm{j}})_{N^{d+1} \times N^{d+1}}, \qquad
A_{\bm{i}, \bm{j}} = \int_{\mathcal{C}_{\mathcal{Y}}} y^{\alpha} \nabla \phi_{\bm{i}}(x,y) \cdot \nabla\phi_{\bm{j}}(x,y) \d x \d y,\]
 \[\bb{b} = (b_{\bm{j}})_{N^{d+1} \times 1}, \qquad
b_{\bm{j}} = d_s \langle f, \text{tr} \phi_{\bm{j}} \rangle = d_s \int_{\Omega} f(x) \phi_{\bm{j}}(x,y)|_{y=0} \,\d x,\]
with $N = M+1$ being the number of points along each direction. Note that the known entries of $ \bb{u} $ corresponding to the Dirichlet boundary conditions should be appropriately adjusted. The approximate solution is obtained by
\[{(\mathfrak{u}_{\mathcal{Y}})}_h: = \sum_{\bb{i}} u_{\mathcal{Y},\bb{i}} \phi_{\bb{i}}, \qquad \bb{u}_{\mathcal{Y}} = (u_{\mathcal{Y},\bb{i}}). \]

To estimate the numerical complexity, we need the following properties of matrix $A$.

\begin{theorem}\label{thm:kappa}
Let $\Omega = (0,1)^d$ and consider the piecewise linear element along each dimension. Then the sparsity $s(A)$ and condition number $\kappa(A)$ satisfy
\[s(A) = 3^{d+1}, \quad
\|A\| = \mathcal{O}((d+1) h^{d-1}), \quad \kappa(A) = \mathcal{O}((d+1)3^{d+1}  h^{-2}).\]
\end{theorem}
\begin{proof}
The sparsity of $A$ is obvious. To estimate the condition number, we only consider $0<s<\frac12$ or $0<\alpha<1$. The case where $s>\frac12$ can be deduced in a similar manner.
For any $\bb{u}$,
\[\bm{u}_{\mathcal{Y}}^\top A \bm{u}_{\mathcal{Y}} =
\int_{\mathcal{C}_{\mathcal{Y}}} y^{\alpha} \nabla {\mathfrak{u}_{\mathcal{Y}}} (x,y) \cdot \nabla {\mathfrak{u}_{\mathcal{Y}}} (x,y) \d x \d y
\le \mathcal{Y}^{\alpha}\int_{\mathcal{C}_{\mathcal{Y}}} \nabla {\mathfrak{u}_{\mathcal{Y}}} (x,y) \cdot \nabla {\mathfrak{u}_{\mathcal{Y}}} (x,y) \d x \d y
= \mathcal{Y}^{\alpha}\bm{u}_{\mathcal{Y}}^\top G_{\mathcal{Y}} \bm{u}_{\mathcal{Y}}, \]
where
\[G_{\mathcal{Y}} = (G^{\mathcal{Y}}_{\bm{i}, \bm{j}})_{N^{d+1} \times N^{d+1}}, \qquad
G^{\mathcal{Y}}_{\bm{i}, \bm{j}} = \int_{\mathcal{C}_{\mathcal{Y}}} \nabla \phi_{\bm{i}}(x,y) \cdot \nabla\phi_{\bm{j}}(x,y) \d x \d y.\]
On the other hand, by the second integral mean-value theorem, there exists $\xi \in [0, \mathcal{Y}]$ such that
\begin{align*}
\bm{u}_{\mathcal{Y}}^\top A \bm{u}_{\mathcal{Y}}
= \int_{\Omega} \Big( \int_0^{\mathcal{Y}} y^{\alpha}  |\nabla {\mathfrak{u}_{\mathcal{Y}}}(x,y)|^2 \d y \Big) \d x
= \int_{\Omega} \Big( {\mathcal{Y}}^{\alpha} \int_{\xi}^{\mathcal{Y}} |\nabla {\mathfrak{u}_{\mathcal{Y}}}(x,y)|^2 \d y \Big) \d x,
\end{align*}
which gives $\bm{u}_{\mathcal{Y}}^\top A \bm{u}_{\mathcal{Y}} = {\mathcal{Y}}^{\alpha} \bm{u}_{\mathcal{Y}}^\top G_{\xi} \bm{u}_{\mathcal{Y}}$, where
\[G_{\xi} = (G^{\xi}_{\bm{i}, \bm{j}})_{N^{d+1} \times N^{d+1}}, \qquad
G^{\xi}_{\bm{i}, \bm{j}} = \int_{\Omega} \int_{\xi}^{\mathcal{Y}} \nabla \phi_{\bm{i}}(x,y) \cdot \nabla\phi_{\bm{j}}(x,y) \d x \d y.\]
The above discussion implies that
\[\kappa(A) \le \frac{\lambda_{\max} (G) }{\lambda_{\min} (G_{\xi})}.\]
As analyzed in \cite[Theorem 3.5]{HJZ2024multiscale} for elliptic PDEs with variable coefficient, we may conclude that
\[\lambda_{\max} (G) \lesssim (d+1) h^{d-1}, \qquad \lambda_{\min} (G_{\xi}) \gtrsim (\frac{h}{3})^{d+1}, \]
which leads to the desired estimate.
\end{proof}

\section{Complexity analysis} \label{sec:complexity}

\subsection{Block encoding of the stiffness matrix} \label{subsec:BE}

In the following, we assume $\Omega = (0,1)^d$. In the $x_k$-direction, we can introduce univariate stiffness and mass matrices:
\begin{equation}\label{AkSk}
A^{(k)} = \Big(\int_0^1 \frac{\partial \varphi_{i_k} (x_k)}{\partial x_k} \frac{\partial \varphi_{j_k} (x_k)}{\partial x_k}\d x_k\Big)_{N\times N}, \qquad
 S^{(k)} = \Big(\int_0^1 \varphi_{i_k} (x_k) \varphi_{j_k} (x_k) \d x_k\Big)_{N\times N}
\end{equation}
 for $k=1,2,\cdots,d$.
In the extended direction, we define weighted stiffness and mass matrices as
\[A^{(d+1)} = \Big(\int_0^{\mathcal{Y}} y^{\alpha}\frac{\partial \psi_{j} (y)}{\partial y} \frac{\partial \psi_{j} (y)}{\partial y}\d y\Big)_{N\times N}, \qquad
 S^{(d+1)} = \Big(\int_0^{\mathcal{Y}} y^{\alpha} \psi_{j} (y) \psi_{j} (y)\d y\Big)_{N\times N}.\]
One can find that the Galerkin stiffness matrix $A$ satisfies
\begin{align*}
A_{\bb{i}, \bb{j}}
& = A_{i_1,j_1} S_{i_2,j_2}\cdots S_{i_{d+1},j_{d+1}} + S_{i_1,j_1} A_{i_2,j_2}\cdots S_{i_{d+1},j_{d+1}} + \cdots
+ S_{i_1,j_1} S_{i_2,j_2}\cdots A_{i_{d+1},j_{d+1}}\\
& =: B_{\bb{i}, \bb{j}} ^{(1)} + B_{\bb{i}, \bb{j}} ^{(2)} + \cdots + B_{\bb{i}, \bb{j}} ^{(d+1)}.
\end{align*}
Noting that
\begin{align*}
B ^{(1)}
& = \sum_{\bb{i}, \bb{j}} B_{\bb{i}, \bb{j}} ^{(1)} \ket{\bb{i}}\ket{\bb{j}}
 = (\sum_{i_1,j_1}A_{i_1,j_1} \ket{i_1,j_1} ) \otimes \cdots \otimes (\sum_{i_{d+1},j_{d+1} }S_{i_{d+1},j_{d+1}} \ket{i_{d+1},j_{d+1}} )\\
& = A^{(1)} \otimes S^{(2)} \otimes \cdots \otimes S^{(d+1)},
\end{align*}
we similarly get the following Kronecker product representation of $A$:
\begin{equation} \label{KronA}
A = A^{(1)} \otimes S^{(2)} \otimes \cdots \otimes S^{(d+1)} + \cdots + S^{(1)} \otimes \cdots \otimes S^{(d)} \otimes A^{(d+1)}.
\end{equation}

\begin{figure}[!htb]
  \centering
  \includegraphics[width=0.75\textwidth]{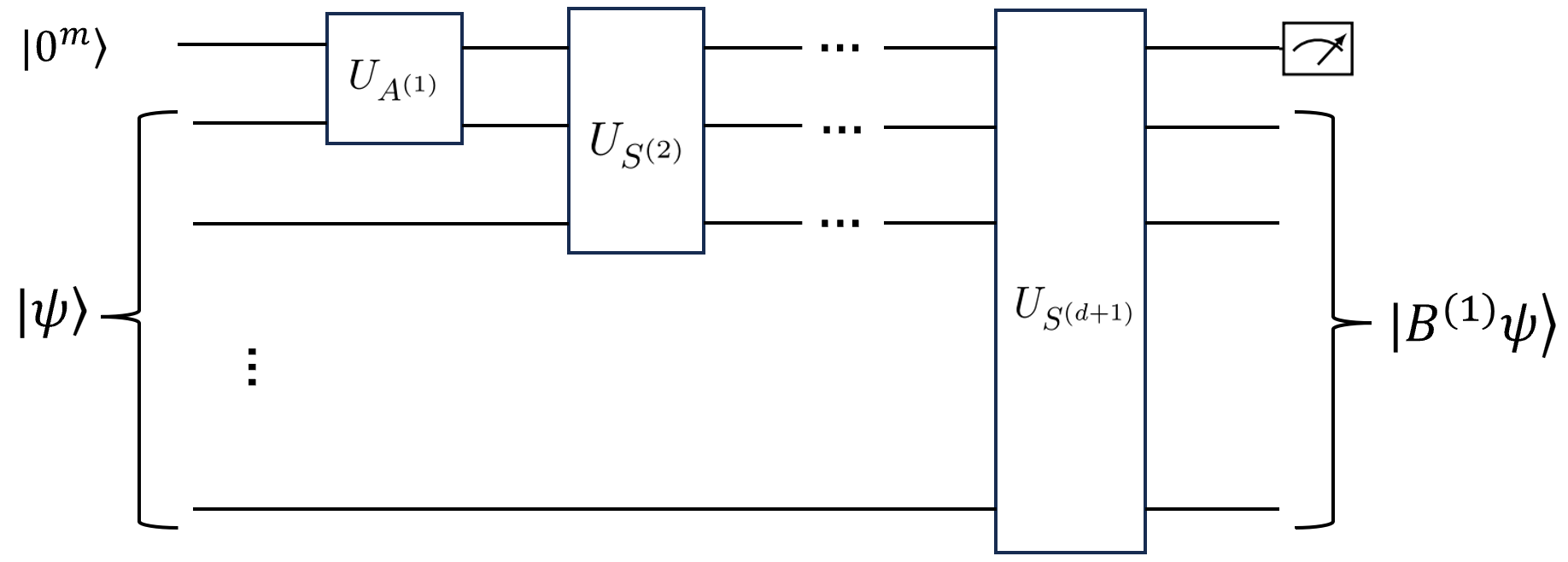}\\
  \caption{Block-encoding of $B^{(1)} = A^{(1)} \otimes S^{(2)} \otimes \cdots \otimes S^{(d+1)}$}\label{fig:BE1}
\end{figure}

Let $ U_{A^{(k)}} $ and $ U_{S^{(k)}} $ denote the block-encodings of $ A^{(k)} $ and $ S^{(k)} $, respectively. We can block-encode $ B^{(1)} = A^{(1)} \otimes S^{(2)} \otimes \cdots \otimes S^{(d+1)} $ as shown in Fig.~\ref{fig:BE1}, where $ U_{S^{(2)}}, \dots, U_{S^{(d+1)}} $ are naturally expanded using existing qubits, while retaining the same ancilla qubits used for $ U_{A^{(1)}} $. The block-encodings of $ B^{(2)}, \dots, B^{(d+1)} $ can be constructed in a similar manner. We can then apply the LCU procedure to encode $ A $.

Although the sparsity $ s(A) = 3^{d+1} $, the complexity of block-encoding $ A $ from the sparse query models for $ A^{(k)} $ and $ S^{(k)} $ has only a quadratic dependence on $ d $.
%As a result, the overall time complexity of the Schr\"odingerization based approach exhibits only a quadratic dependence on $d$ if we apply the established block-encoding input model.

\subsection{Trace projection}

Let $0=y_0<y_1<\cdots<y_{N-1}= \mathcal{Y}$ be the grid points in the extended dimension. The corresponding basis functions  are denoted by $\psi_0(y), \cdots, \psi_{N-1}(y)$, respectively. Noting that
\[\psi_0(0)=1, \quad \psi_1(0)= \cdots = \psi_{N-1}(0),\]
we can express the FEM solution for the original equation as
\[u_h = \text{tr}{(\mathfrak{u}_{\mathcal{Y}})}_h: = \sum_{i_1,\cdots,i_d=0}^{N-1}  u_{i_1,\cdots,i_d, i_{d+1}} \phi_{i_1,\cdots,i_d, i_{d+1}}, \qquad i_{d+1} = 0. \]
Let
\begin{equation}\label{uhFEM}
u_h = \sum_{i_1,\cdots,i_d=0}^{N-1}  u_{i_1,\cdots,i_d} \phi_{i_1,\cdots,i_d}.
\end{equation}
One may get $u_{i_1,\cdots,i_d} = u_{i_1,\cdots,i_d, i_{d+1}=0}$.
This means
\begin{equation}\label{uhvectorFEM}
\bb{u}_h = \sum_{i_1,\cdots,i_d}^{N-1}  u_{i_1,\cdots,i_d} \ket{i_1,\cdots,i_d} = \sum_{i_1,\cdots,i_d}^{N-1}  u_{i_1,\cdots,i_d, i_{d+1}=0} \ket{i_1,\cdots,i_d, i_{d+1}=0},
\end{equation}
where $\bb{u}_h$ denotes the solution vector of $u_h$ in the basis of the $x$-domain, which can be obtained from $\bb{u}_{\mathcal{Y}}$ by projecting onto $\ket{0}$ in the extended dimension. The probability for obtaining $\bb{u}_h$ is $\frac{\|\bb{u}_h\|^2}{\|\bb{u}_{\mathcal{Y}}\|^2}$. The success probability can be raised to $\Omega(1)$ by using $\mathcal{O}(g)$ rounds of amplitude amplification, where $g = \mathcal{O}(\frac{\|\bb{u}_{\mathcal{Y}}\|}{\|\bb{u}_h\|})$.

For the upper bound of $g$, one can apply the standard norm equivalence in FEMs.
\begin{lemma}
Under the condition of Theorem \ref{thm:kappa}, there holds
\begin{align*}
& (\frac{h}{3})^{-\frac{d}{2}} \|u_h\|_{L^2} \le  \|\bb{u}_h\| \le  h^{-\frac{d}{2} } \|u_h\|_{L^2}, \\
& (\frac{h}{3})^{-\frac{d+1}{2}} \|(u_{\mathcal{Y}})_h\|_{L^2} \lesssim  \|\bb{u}_{\mathcal{Y}}\| \lesssim  h^{-\frac{d+1}{2} } \|(u_{\mathcal{Y}})_h\|_{L^2},
\end{align*}
implying
\begin{equation}\label{grepeat}
g = \mathcal{O}( 3^{\frac{d}{2}} h^{-\frac{1}{2}}).
\end{equation}
\end{lemma}
\begin{proof}
Let $S$ be the mass matrix obtained from the basis functions in $\Omega$. From \eqref{uhFEM} and following the derivation in \eqref{KronA}, we get
\[
\|u_h\|_{L^2}^2 = \bb{u}_h^T S \bb{u}_h, \qquad S = S^{(1)} \otimes S^{(2)} \otimes \cdots \otimes S^{(d)},
\]
where $S^{(k)}$ is defined in \eqref{AkSk} for $k=1,2,\cdots,d$. According to the computation in \cite[Theorem 3.5]{HJZ2024multiscale}, the eigenvalues satisfy
\begin{equation}\label{Sklambda}
\frac{h}{3} \le \lambda(S^{(k)}) \le h, \qquad k = 1,2,\cdots, d,
\end{equation}
leading to
\[(\frac{h}{3})^d \le \lambda(S) \le h^d.\]
Therefore,
\[(\frac{h}{3})^d \|\bb{u}_h\|^2 \le \|u_h\|_{L^2}^2 \le h^d \|\bb{u}_h\|^2.\]

For the second inequality, one can apply the second integral mean-value theorem as done in Theorem \ref{thm:kappa} and obtain
\[\frac{h}{3} \lesssim \lambda(S^{(d+1)}) \lesssim h.\]
This completes the proof.
\end{proof}

\subsection{Time complexity}

The linear system $A\bb{x} = \bb{b}$ for the FEM disretization of the extension problem is given in \eqref{AxbFEM}, where
$\bb{x} = \bb{u}_{\mathcal{Y}}$. For our quantum algorithm, the solution to the linear system is interpreted as the steady-state solution to the ODEs in \eqref{reformulationODE}. The corresponding solution vector is denoted as $\bb{x}_T = \bb{u}_{\mathcal{Y}}^t$, which is then approximated by the solution vector $\bb{x}_T^{\text{d}} = \bb{u}_{\mathcal{Y}}^q$ in the Schr\"odingerization procedure. The vectors in the $x$-direction are denoted by $\bb{u}_h$, $\bb{u}_h^t$ and $\bb{u}_h^q$, respectively.
We can accordingly define the approximations of $u$:
\begin{equation}\label{functionexpression}
u_h = \sum_{\bm{i}_x } u_{h,\bm{i}_x} \phi_{\bm{i}_x}(x), \qquad
 %u_h^t = \sum_{\bm{i}_x } u_{h,\bm{i}_x}^t \phi_{\bm{i}_x}(x),\quad
 u_h^q = \sum_{\bm{i}_x } u_{h,\bm{i}_x}^q \phi_{\bm{i}_x}(x),
\end{equation}
where $\bb{i}_x = (i_1,\cdots,i_d)$ and $\phi_{\bm{i}_x}(x)$ is a basis function in the $x$-direction.
Here, we omit the definition of $ u_h^t $ as our error analysis does not involve it.

The goal of this section is to determine the computational cost of obtaining a quantum solution $ u^q $ that satisfies the following condition
\[\|u - u_h^q\|_{L^2(\Omega)} \le \delta = \delta(h).\]
The error can be divided into two parts:
\[ \|u - u_h^q\|_{L^2}  \le  \|u - u_h\|_{L^2}
+ \|u_h - u_h^q\|_{L^2}.\]
We will require each part of the error to be $\mathcal{O}(\delta)$.

The first term, which corresponds to the finite element error, can be estimated using Lemma \ref{lem:err}:
\[\|u - u_h\|_{L^2(\Omega)}  \le \|u - u_h\|_{\mathbb{H}^s(\Omega)} \lesssim  | d \log (h^{-1} )| ^s h =:\delta. \]
For the second term, by the norm equivalence,
\[\|u_h - u_h^q\|_{L^2(\Omega)}
\lesssim h^{\frac{d}{2}} \|\bb{u}_h - \bb{u}_h^q\|
\le h^{\frac{d}{2}} \|\bb{u}_{\mathcal{Y}} - \bb{u}_{\mathcal{Y}}^q\| = h^{\frac{d}{2}} \|\bb{x}  - \bb{x}_T^{\text{d}}\|.\]
Requiring the right-hand side to be $\mathcal{O}(\delta)$ and noting that
\begin{align*}
\|\ket{\bb{x}} - \ket{\bb{x}_T^{\text{d}}}\|
 \le 2 \frac{\|\bb{x} - \bb{x}_T^{\text{d}}\|}{\|\bb{x}\|}
 \le 2 \frac{\|\bb{x} - \bb{x}_T^{\text{d}}\|}{\|\bb{u}_{\mathcal{Y}}\|}
 \le 2 \frac{\|\bb{x} - \bb{x}_T^{\text{d}}\|}{\|\bb{u}_h\|},
\end{align*}
we can take
\[\varepsilon \lesssim \frac{\delta h^{-\frac{d}{2}}}{\|\bb{u}_h\|}
\lesssim \frac{ \delta }{ \|u\|_{L^2} } = \mathcal{O}( \delta ) \]
in Theorem \ref{thm:SchrCost}.

\begin{theorem}\label{thm1:SchrCost}
Let $ \varepsilon = \mathcal{O}( | d \log (h^{-1} )|^s h ) $. Under the conditions of Theorem \ref{thm:SchrCost}, there exists a quantum algorithm that prepares an $ \varepsilon $-approximation of the state $ \ket{\bb{x}} $, denoted by $ \ket{\bb{x}_T^{\text{d}}} $, with $ \Omega(1) $ success probability and a flag indicating success, using $ \widetilde{\mathcal{O}}( d 3^{\frac{3d}{2}} h^{-2.5} ) $ queries to the block-encoding oracle of $ A $ and the state preparation oracle for $ \bb{b} $. The unnormalized vector $ \bb{x}_T^{\text{d}} $ approximates $ \bb{u}_{\mathcal{Y}} $, and extracting in the $ x $-direction yields an approximation of $ \bb{u}_h $, denoted as $ \bb{u}_h^q $, which satisfies $ \| u - u_h^q \|_{L^2(\Omega)} \lesssim \varepsilon $, with $ u_h^q $ defined in \eqref{functionexpression}.
\end{theorem}

\begin{proof}
The result is obtained by combining Theorem \ref{thm:SchrCost}, the VTAA procedure (see Remark \ref{rem:VTAA}), Theorem \ref{thm:kappa} and Eq.~\eqref{grepeat}.
\end{proof}

On the other hand, the time complexity of the classical CG  method for solving \eqref{AxbFEM} is
\[C_{\text{CG}} = \mathcal{O}(N^{d+1} s \sqrt{\kappa} \log (1/\varepsilon)) =
\widetilde{\mathcal{O}}(d^{\frac12} 3^{\frac{3d}{2}} h^{-d-2}).\]
Considering that VTAA requires substantial modifications to the LCU or quantum signal processing algorithms, which introduces significant overhead itself, we can derive the query complexity, neglecting the VTAA procedure, as
\[C_{\text{Quantum}} = \widetilde{\mathcal{O}}(d^2 3^{\frac{5d}{2}} h^{-4.5}).\]
The quantum algorithm exhibits a dependence on \( h^{-1} \) that is {\it independent} of the dimension, whereas the classical method has an exponential dependence, particularly in high dimensions. This demonstrates the exponential quantum advantage for high-dimensional fractional problems, where classical methods suffer from the curse of dimensionality.

\section{Numerical simulation} \label{sec:numerical}

This section reports the performance of our proposed method and compares it with the method in \cite{Acosta2015Fractional}, which is based on a nonlocal variational formulation.

\subsection{The Caffarelli-Silvestre extension}

\subsubsection{The finite element discretization}

Let $\Omega = (-1,1)^d$. The exact solution and the right-hand side are given by
\[u(x_1,\cdots,x_d) = \sin (\pi x_1) \cdots \sin (\pi x_d), \qquad f(x_1,\cdots,x_d) = (d\pi^2)^s u(x_1,\cdots,x_d).\]
Unless otherwise specified, we take $N_p = 2^{11}$ for the discrete Fourier transform in the implementation of the Schr\"odingerization method. The evolution time is selected as $T = 15$, which is sufficient for the steady-state solution. For the truncation in the extended dimension, we take $\mathcal{Y} = 10$.

%\subsection{One dimensional case}

\begin{figure}[!htb]
		\centering
        \subfigure[Extended solution $\mathfrak{u}(x,y)$]{\includegraphics[scale=0.45]{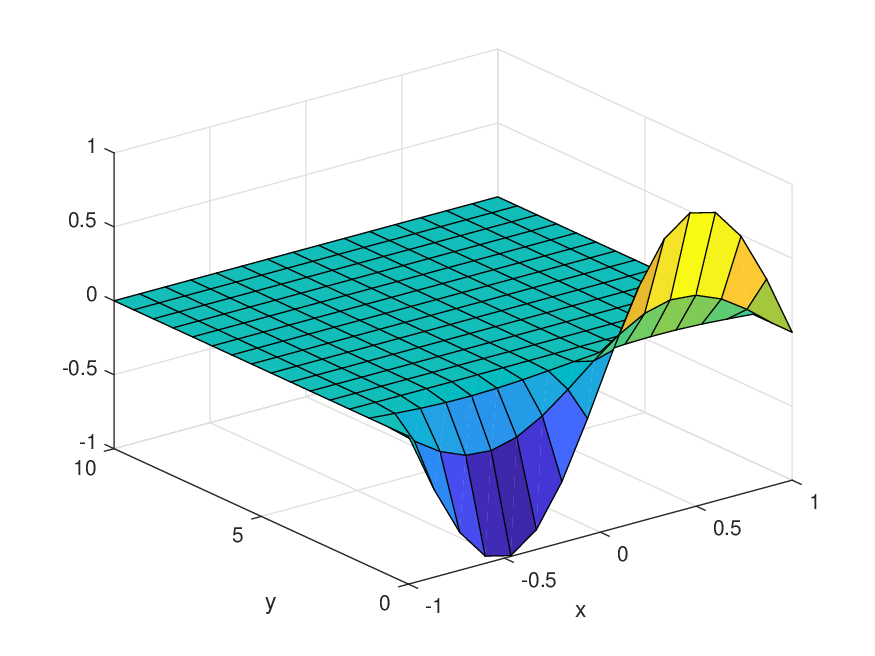}}
		\subfigure[Solution $u(x)$]{\includegraphics[scale=0.45]{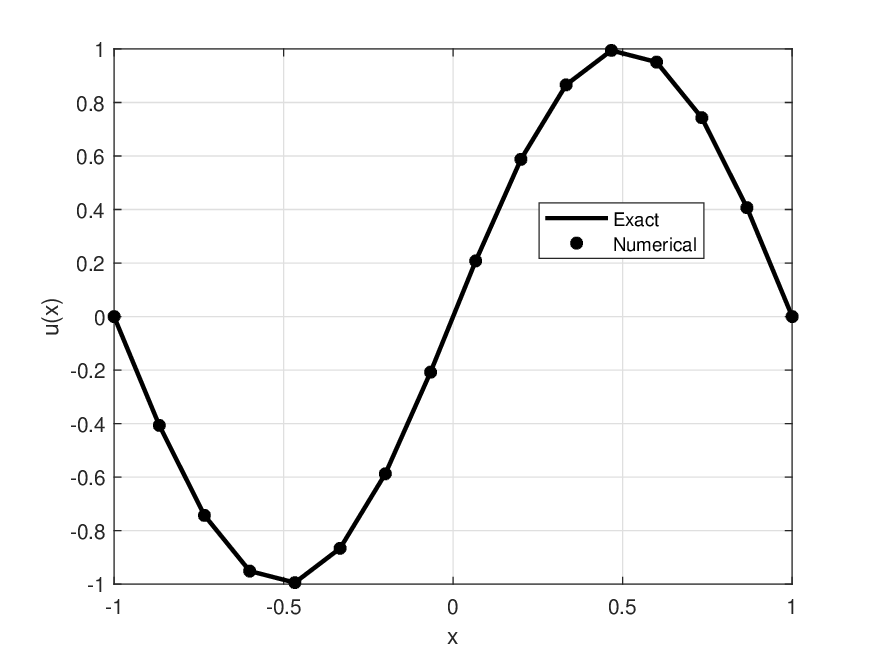}}\\
		\caption{Numerical solutions in 1D with uniform partition ($s=0.2$)}\label{fig:CS1d}
	\end{figure}

We begin by considering the one-dimensional case. The numerical results, obtained with a uniform partition of the interval $[0, \mathcal{Y}]$, are presented in Fig.~\ref{fig:CS1d} for $s=0.2$. The extended solution $\mathfrak{u}(x,y)$ is shown in Fig.~\ref{fig:CS1d}a, where we observe that the solution in the extended dimension decays rapidly to zero. The original solution $u(x) = \mathfrak{u}(x,y)|_{y=0}$ is also displayed in Fig.~\ref{fig:CS1d}, where we see that the Schr\"odingerization based approach provides an accurate approximation.

\begin{figure}[!htb]
		\centering
        \subfigure[Extended solution $\mathfrak{u}(x,y)$]{\includegraphics[scale=0.45]{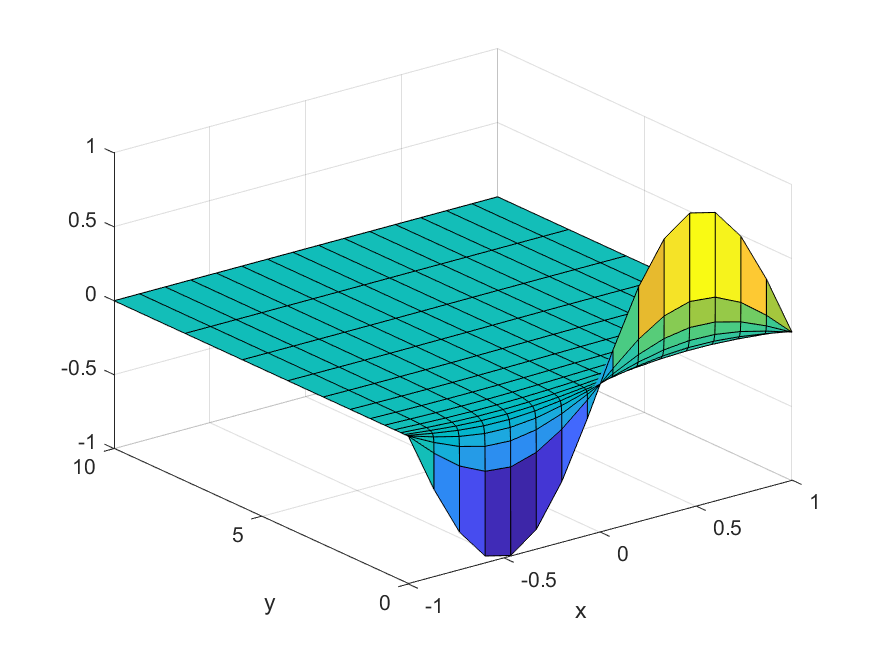}}
		\subfigure[Solution $u(x)$]{\includegraphics[scale=0.45]{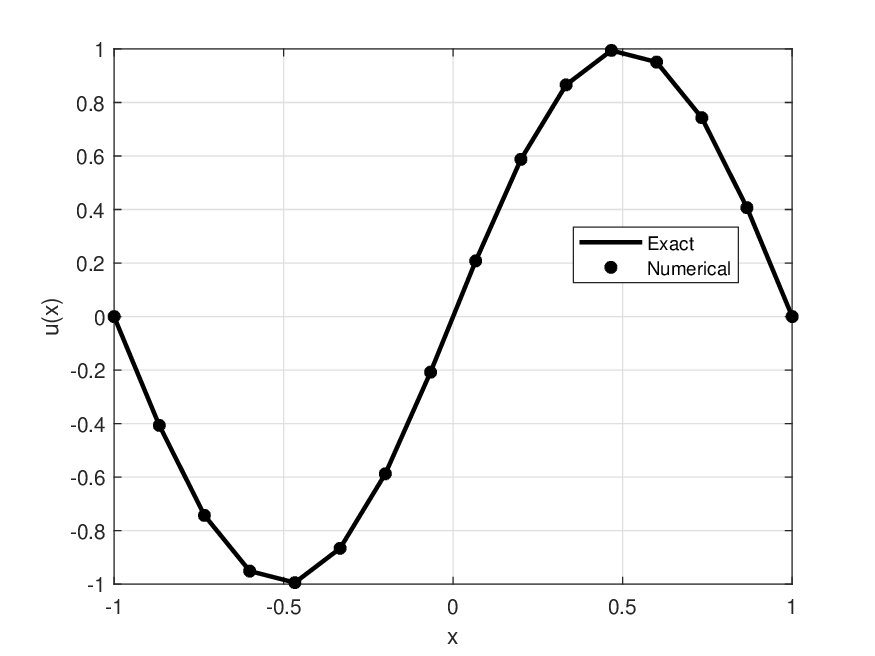}}\\
		\caption{Numerical solutions in 1D with graded partition ($s=0.8$)}\label{fig:CS1dgrad}
	\end{figure}

We repeat the test with a graded partition as described in \eqref{graded}, using $\gamma = 4$ for $s = 0.8$. Similar results are observed and displayed in Fig.~\ref{fig:CS1dgrad}.

%\subsection{Two dimensional case}

\begin{figure}[!htb]
		\centering
        \subfigure[Exact]{\includegraphics[scale=0.45]{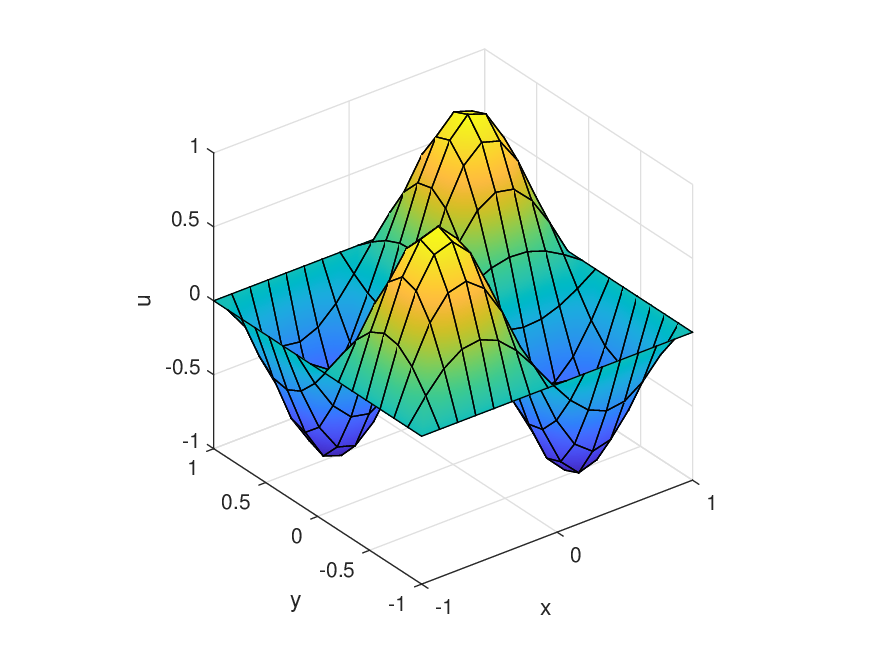}}
		\subfigure[Schr\"odingerization]{\includegraphics[scale=0.45]{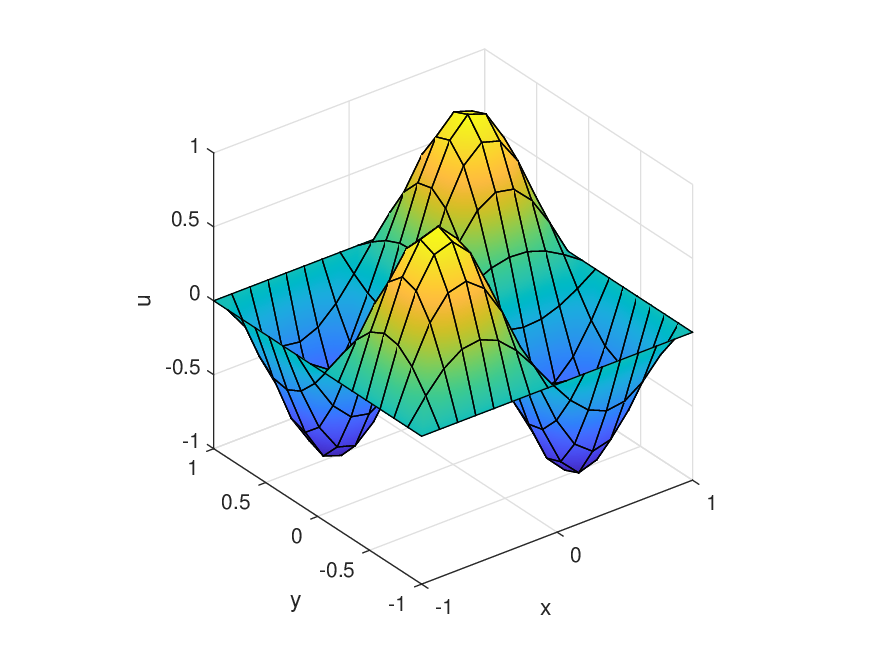}}\\
		\caption{Numerical solutions of $u$ based on the Caffarelli-Silvestre extension ($s=0.2$)}\label{fig:CS}
	\end{figure}

\begin{figure}[!htb]
		\centering
        \includegraphics[scale=0.45]{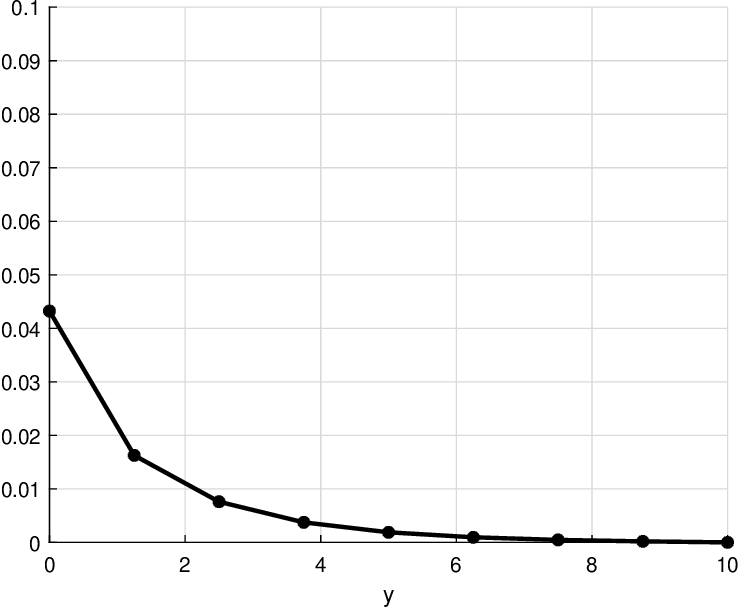}\\
		\caption{Solution profile along the extended direction at fixed $(x_1,x_2) = (0,0)$}\label{fig:CSuz}
	\end{figure}

For the test in two dimensions, we present the exact and numerical solutions for $u(x_1,x_2) = \mathfrak{u}(x_1,x_2,y)|_{y=0}$ with $s = 0.2$ in Fig.~\ref{fig:CS}, where we observe that the Schr\"odingerization-based approach provides an accurate approximation. Additionally, we display the solution profile of $\mathfrak{u}(x_1,x_2,y)$ at the fixed point $(x_1,x_2) = (0,0)$ in Fig.~\ref{fig:CSuz}. As expected, the solution along the $y$-direction decays to zero.

\subsubsection{The finite difference discretization}

\begin{figure}[!htb]
		\centering
        \includegraphics[scale=0.45]{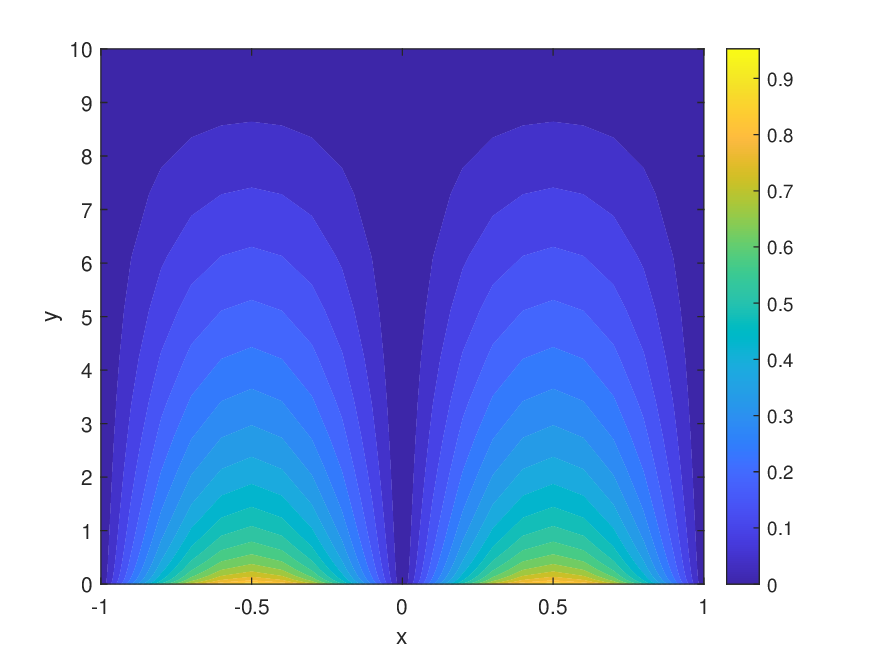}\\
		\caption{Solution profile along the extended direction at fixed $(x_1,x_2) = (0,0)$ with $s = 0.8$}\label{fig:CSFDM}
	\end{figure}

The Schr\"odingerization approach is also applicable to finite difference discretizations for the extension problem. For this test, we use central differences for the second-order derivatives, with the contour map of $|\mathfrak{u}(x,y)|$ in 1D shown in Fig.~\ref{fig:CSFDM}. It is evident that the solution decays rapidly to zero along the $y$ direction.

\subsection{Nonlocal variational formulation} \label{subsec:nonlocal}

Ref.~\cite{Acosta2015Fractional} deals with the integral version of the Dirichlet homogeneous fractional Laplace equation, and proposes a linear finite element method. The nonlocal variational formulation reads: Find $u \in \widetilde{H}^s(\Omega)$ such that
\[\frac{C_{d,s}}{2} \iint_{\mathbb{R}^d \times \mathbb{R}^d}
\frac{(u(x)-u(y))(v(x)-v(y))}{|x-y|^{d + 2s}} \d x \d y = \int_{\Omega} f v \d x, \qquad v \in \widetilde{H}^s(\Omega), \]
where $\widetilde{H}^s(\Omega) = \{ v \in H^s(\mathbb{R}^d): \text{supp}~v \subset \bar{\Omega}\}$.

\begin{figure}[!htb]
		\centering
		\subfigure[nonlocal variational formulation]{\includegraphics[width=0.3\textwidth,trim=40 0 40 0,clip]{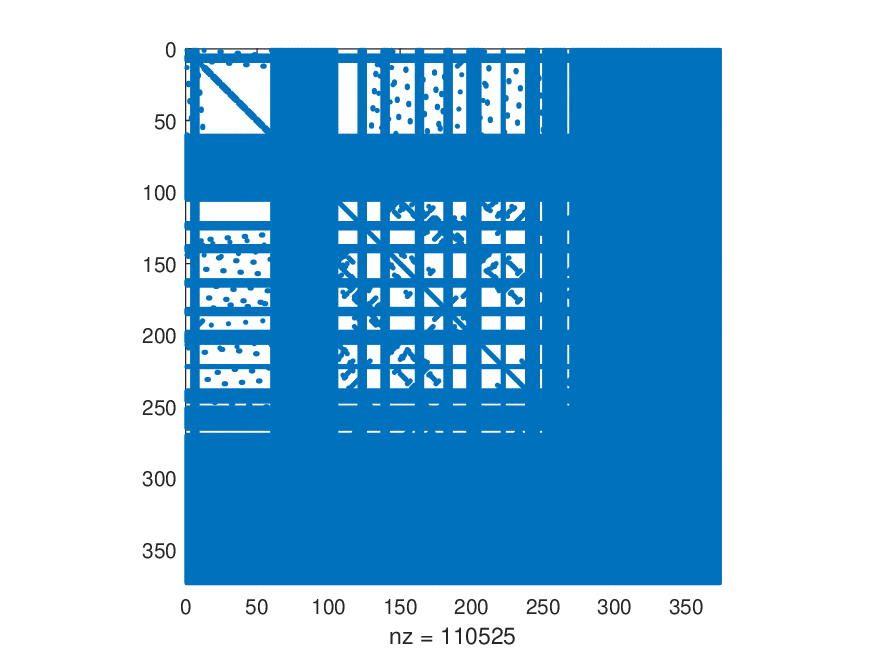}}
        \subfigure[Extension in 1D]{\includegraphics[width=0.3\textwidth,trim=40 0 40 0,clip]{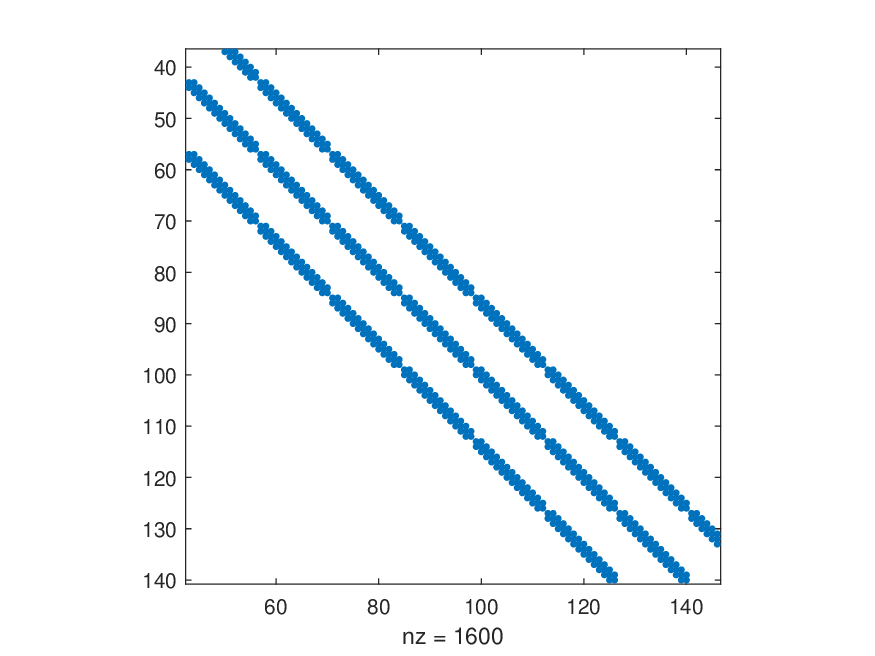}}
        \subfigure[Extension in 2D (zoomed)]{\includegraphics[width=0.3\textwidth,trim=40 0 40 0,clip]{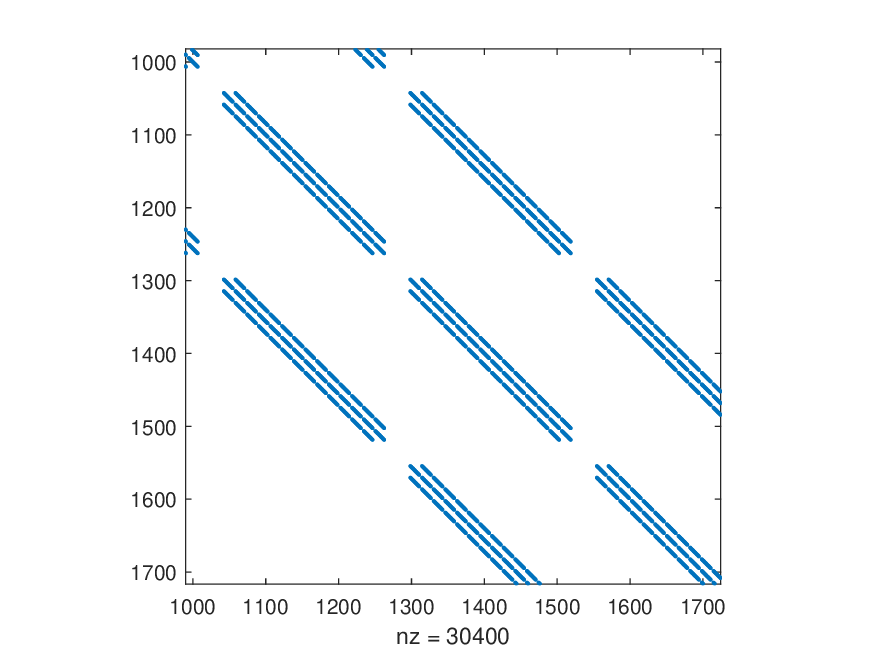}}\\
		\caption{Sparsity pattern for linear systems}\label{fig:sparsityNonlocal}
	\end{figure}

We continue to consider the test in two dimensions. As discussed in \cite{Acosta2017codeFractional}, when implementing the finite element method, it is helpful to introduce a ball $ B $ that contains $ \Omega $ when dealing with integrals over $ \Omega^c = \mathbb{R}^2 \backslash \Omega $, especially when $ \Omega $ is not a ball. In this case, the nodes in the auxiliary domain $ B \backslash \Omega $ are treated as Dirichlet nodes. Since the bilinear form involves integrals in $ \mathbb{R}^d $, its Galerkin discretization results in dense matrices, as shown in Fig.~\ref{fig:sparsityNonlocal}a. This introduces challenges in preparing the input model for the linear system, especially when $\Omega$ is a general bounded domain. The sparsity pattern for the Caffarelli-Silvestre extension differs significantly from this nonlocal form, as illustrated in Figs.~\ref{fig:sparsityNonlocal}b and c, where $s(A) = 3^2$ and $s(A) = 3^3$ are observed, consistent with the result in Theorem \ref{thm:kappa}. Note that in the sparsity plot, each diagonal line corresponds to a sparsity of 3.

\begin{figure}[!htb]
		\centering
        \subfigure[Direct solver]{\includegraphics[scale=0.45]{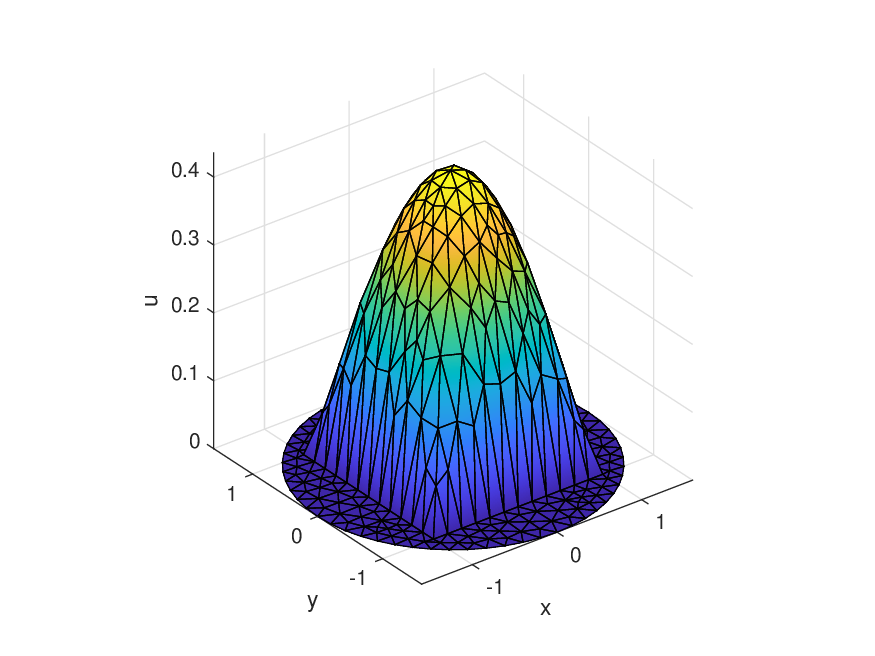}}
		\subfigure[Schr\"odingerization]{\includegraphics[scale=0.45]{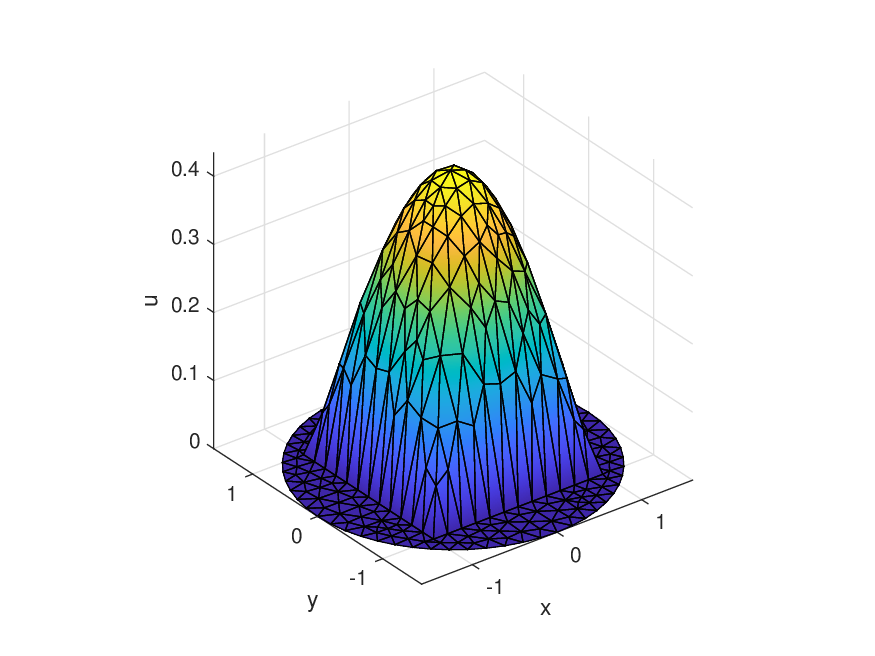}}\\
		\caption{Numerical solutions of $u$ based on nonlocal variational formulation ($s=0.8$)}\label{fig:solNonlocal}
	\end{figure}

Let $ \Omega = (-1,1)^2 $. The ball $ B $ used in the numerical implementation is centered at $ (0,0) $ with a radius of 1.5. The results for $ s = 0.8 $ are presented in Fig.~\ref{fig:solNonlocal}, with both direct and Schr\"odingerization solvers applied. Since the enlarged domain is a circle, a linear finite element method on a triangulation is used for convenience. It is observed that the Schr\"odingerization-based approach produces accurate approximate results.

\begin{remark} \label{rem:fast}
Ref.~\cite{Yang2023fastNonlocal} proposes a fast Q1 finite element method based on a weighted trapezoidal rule for the nonlocal variational formulation. The resulting matrix is a symmetric block Toeplitz matrix, which can be efficiently computed using FFT and its inverse transform. This property was also observed in \cite{Hao2021fractionalFDM,YangChen2023fastNonlocalFDM} for finite difference discretizations in multiple dimensions. Consequently, we may provide an efficient method for block-encoding the linear system of the nonlocal form in some cases, even though the stiffness matrix is dense.
\end{remark}

\section{Conclusion}

In this article we developed a quantum algorithm for solving the Dirichlet homogeneous fractional Poisson equation.
The algorithm combines two key techniques: (1) the Caffarelli-Silvestre extension, which reduces the nonlocal fractional Laplacian to a local elliptic problem in one higher dimension, and (2) the Schr\"odingerization method, which converts the resulting PDE with finite element discretizations into a Hamiltonian simulation problem.
We conducted a thorough analysis of the time complexity and demonstrated that exponential quantum advantages can be achieved with respect to the
inverse of the mesh size in high dimensions.
Specifically, while the classical method requires \(\widetilde{\mathcal{O}}(d^{1/2} 3^{3d/2} h^{-d-2})\) operations, the quantum counterpart requires \(\widetilde{\mathcal{O}}(d 3^{3d/2} h^{-2.5})\) queries to the block-encoding input models, with the quantum complexity being independent of the dimension \(d\) in terms of \(h^{-1}\).

\section*{Acknowledgments}
SJ and NL were supported by NSFC grant No. 12341104,
the Shanghai Jiao Tong University 2030 Initiative and the Fundamental Research Funds for the Central Universities. SJ was also partially supported by the NSFC grant No. 12426637 and  the Shanghai Municipal Science and Technology Major Project (2021SHZDZX0102). NL also acknowledges funding from the Science and Technology Program of Shanghai, China (21JC1402900), the Science and Technology Commission of Shanghai Municipality (STCSM) grant no. 24LZ1401200 (21JC1402900) and NSFC grant No.12471411.
YY was supported by NSFC grant No. 12301561, the Key Project of Scientific Research Project of Hunan Provincial Department of Education (No. 24A0100) and the 111 Project (No. D23017).

\bibliographystyle{plain} %plain, unsrt, alpha
\bibliography{Refs}

\end{document}